\documentclass[onefignum,onetabnum]{siamart190516}

\usepackage[normalem]{ulem} 

\usepackage{amssymb}
\usepackage{adjustbox}
\newcommand{\vvvertiii}[1]{{\left\vert\kern-0.15ex\left\vert\kern-0.15ex\left\vert #1 \right\vert\kern-0.15ex\right\vert\kern-0.15ex\right\vert}}

\usepackage{bm}
\usepackage{framed}
\usepackage{mathrsfs}
\usepackage{stmaryrd}

\usepackage{tikz}
\usepackage{environ}
\NewEnviron{mybox}[1][]{%
 \noindent
 \begin{tikzpicture}
   \node[minimum width=\linewidth-0.4pt,draw,text width=\linewidth-12pt] (a){\BODY};
   \node[fill=white,yshift=0.5ex] at (a.north) {#1};
 \end{tikzpicture}}

\usepackage[english]{babel}

\usepackage{soul}

\definecolor{Cerulean}{HTML}{007BA7}
\definecolor{Orchid}{HTML}{BA55D3}
\definecolor{ForestGreen}{HTML}{228B22}
\definecolor{BlueViolet}{HTML}{473992}

\newtheorem{assumption}{Assumption}
\newtheorem{remark}{Remark}

\newcommand{\SM}{H} 
\newcommand{\VM}{P} 
\newcommand{\B}{B} 
\newcommand{\KLP}{K} 
\newcommand{\CSconstant}{\gamma}

\newcommand{\uDim}{n} 
\newcommand{\sDim}{d} 
\newcommand{\thickness}{\eta}
\newcommand{\moment}{B}
\newcommand{\ersatz}{\textup{T}}
\newcommand{\bendStrain}{\beta}
\newcommand{\midsurfRot}{\theta}
\newcommand{\bdyTangent}{t}

\newcommand{\applied}[1]{\hat{#1}}
\newcommand{\cTrace}[1]{C_{\textup{tr} #1}}

\newcommand{\cPen}[1]{C_{\textup{pen} #1}}
\newcommand{\cornerSet}[1]{\chi_{#1}}

\newcommand{\physDomain}{\Omega}
\newcommand{\physBoundary}{\Gamma}

\newcommand{\forcing}{\textup{f}}

\newcommand{\linFunctional}{f}

\newcommand{\N}{\mathbb N}
\newcommand{\R}{\mathbb R}

\numberwithin{equation}{section}

\headers{Constructing Nitsche's method for variational problems}{J. Benzaken, J. A. Evans, and R. Tamstorf}

\title{Constructing Nitsche's method for variational problems}

\author{Joseph Benzaken\thanks{Walt Disney Animation Studios, Burbank, CA (\email{joseph.benzaken@disneyanimation.com}, \newline\email{Rasmus.Tamstorf@disneyanimation.com}).}
\and John A. Evans\thanks{Ann and H.J. Smead Department of Aerospace Engineering Sciences, University of Colorado at Boulder, Boulder, CO (\email{john.a.evans@colorado.edu}).}
\and Rasmus Tamstorf\footnotemark[1]
}

\begin{document}

\maketitle

\begin{abstract}
  Nitsche's method is a well-established approach for weak enforcement of boundary conditions for partial differential equations (PDEs). It has many desirable properties, including the preservation of variational consistency and the fact that it yields symmetric, positive-definite discrete linear systems that are not overly ill-conditioned. In recent years, the method has gained in popularity in a number of areas, including isogeometric analysis, immersed methods, and contact mechanics. However, arriving at a formulation based on Nitsche's method can be a mathematically arduous process, especially for high-order PDEs. Fortunately, the derivation is conceptually straightforward in the context of variational problems. To facilitate the process, we devised an abstract framework for constructing Nitsche's method for these types of problems in \cite{Benzaken2020}. The goal of this paper is to elucidate the process through a sequence of didactic examples. First, we show the derivation of Nitsche's method for Poisson's equation to gain an intuition for the various steps. Next, we present the abstract framework and then revisit the derivation for Poisson's equation to use the framework and add mathematical rigor. In the process, we extend our derivation to cover the vector-valued setting. Armed with a basic recipe, we then show how to handle a higher-order problem by considering the vector-valued biharmonic equation and the linearized Kirchhoff-Love plate. In the end, the hope is that the reader will be able to apply Nitsche's method to any problem that arises from variational principles.
\end{abstract}



\tableofcontents

\section{Introduction}
\label{sec:introduction}

Since the broad adoption of the finite element method in the mid-twentieth century, several different approaches for the enforcement of Dirichlet boundary conditions have been proposed. The most common approach is strong enforcement within the finite element space via interpolation or projection, leading to a standard Bubnov-Galerkin method \cite{hughes2012finite}.  For simple variational problems such as heat conduction and linear elasticity, strong enforcement of boundary conditions is straightforward and leads to a stable and convergent finite element method.  However, for fourth-order problems this is much more difficult since both functional and derivative boundary conditions must be enforced.  For interface problems, it can even lead to a detrimental loss in convergence rate, a phenomenon known as ``boundary locking'' \cite{lew2008discontinuous}. A second common approach to enforcing Dirichlet boundary conditions is to do so weakly with a set of Lagrange multipliers.  The method of Lagrange multipliers is particularly popular in the context of interface problems, and both mortar finite element methods \cite{belgacem1999mortar,wohlmuth2000mortar} and finite element tearing and interconnecting (FETI) methods \cite{farhat2001feti,farhat1991method} involve the introduction of Lagrange multiplier fields.  The primary disadvantage of the Lagrange multiplier method is that it leads to discrete saddle-point problems, and the stability of the method can only be ensured if the approximation spaces for the primal and Lagrange multiplier fields satisfy the Babu\v{s}ka-Brezzi inf-sup condition \cite{brezzi2012mixed}.  A third approach for weak Dirichlet boundary condition enforcement is through the penalty method \cite{Babuska1973}. In this approach, the boundary conditions are enforced via boundary regularization terms with tunable penalty parameters that specify the degree to which the boundary conditions are enforced. However, the penalty method is generally inaccurate unless the penalty parameters are chosen to be very large, and this in turn results in severe ill-conditioning.

Nitsche's method is an alternative technique for the weak enforcement of boundary conditions that yields a formulation that is both consistent and stable, and it also provides optimal convergence rates after discretization. In contrast to the method of Lagrange multipliers, Nitsche's method does not result in a discrete saddle-point problem, and therefore is not subject to the Babu\v{s}ka-Brezzi inf-sup condition. Furthermore, unlike the penalty method, it achieves accurate results using moderately sized penalty parameters, guided by trace inequalities \cite{Evans2013,warburton2003constants}, leading to relatively well-conditioned symmetric positive-definite linear systems after discretization of self-adjoint elliptic PDEs.

\begin{figure}
\centering
    \includegraphics{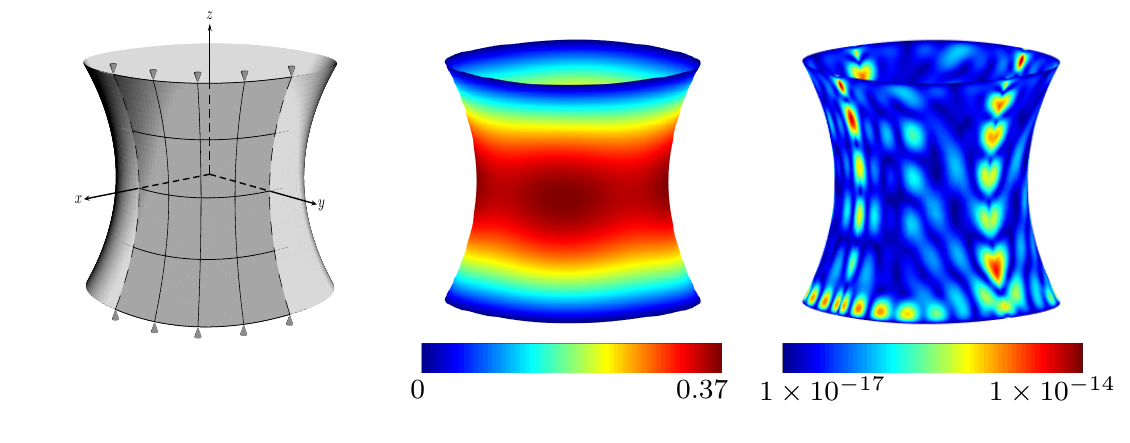}
  \caption{A $p=6$, $4 \times 4$-element B-spline approximation of a thin shell undergoing deformation from a hyperbolic paraboloid to a cylinder. The problem setup (left) shows that the top and bottom are simply supported, and that only one quarter of the shell is considered due to the symmetry of the problem configuration. The symmetry is accounted for through the boundary conditions, and all boundary conditions are enforced using Nitsche's method. The resulting displacement field is plotted over the undeformed geometry (center) with the error between the analytical displacement field and the numerical solution (right). Note that the magnitude of the error is close to machine precision everywhere.}
 \label{fig:shell_ex}
\end{figure}

Nitsche's method was first proposed in 1971 for the weak enforcement of boundary conditions in the finite element approximation of the Poisson problem \cite{Nitsche1971}, but it did not grow in popularity until rather recently with the emergence of isogeometric \cite{embar2010imposing,Apostolatos2014,nguyen2014nitsche,ruess2014weak} and immersed \cite{annavarapu2012robust,hansbo2002unfitted,FernandezMendez2004,kamensky2015immersogeometric,ruess2013weakly,schillinger2012isogeometric} finite element methods.  Isogeometric finite element methods employ spline basis functions originating in computer aided design instead of standard finite element approximation functions.  As a result, these methods offer the possibility of more tightly integrating computer aided design with analysis.  Moreover, as spline basis functions typically exhibit $C^1$ or greater continuity, isogeometric finite element methods can also be applied to the numerical solution of higher-order partial differential equations, such as those governing Kirchhoff-Love shells \cite{kiendl2009isogeometric}.

\begin{figure}
\centering
    \includegraphics[width=5in]{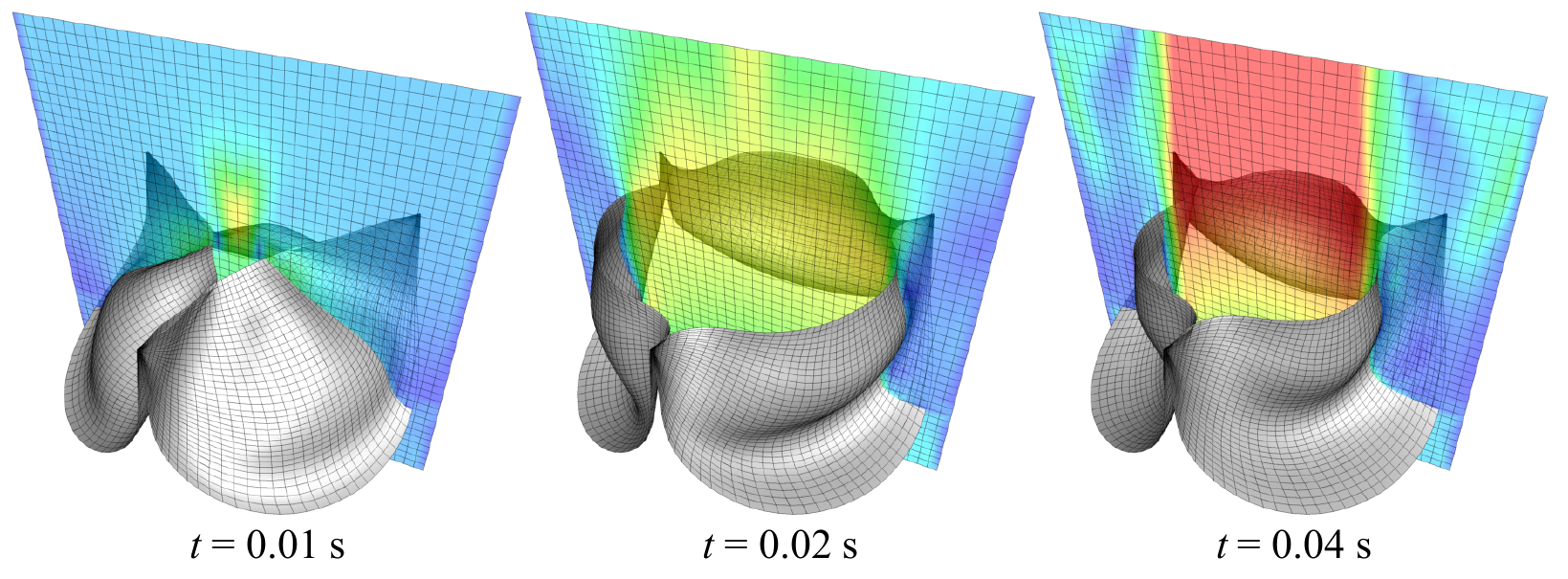}
  \caption{Snapshots of a heart valve opening within a cylindrical fluid domain.  For three representative times, velocity magnitude is plotted on a slice, using a color scale ranging from 0 cm/s (blue) to greater than 200 cm/s (red).  The opening process is simulated using a particular immersed finite element technique known as immersogeometric analysis.  Valve deformation is modeled using Kirchhoff-Love theory, while blood flow is modeled using the incompressible Navier-Stokes equations.  The valve and cylindrical fluid domains are discretized using NURBS and divergence-conforming B-splines, respectively.  Continuity of velocity across the fluid-structure interface is weakly enforced using Nitsche's method.  This figure is reprinted with permission from \cite{kamensky2017immersogeometric}.}
 \label{fig:valve_ex}
\end{figure}

In the case of immersed finite element methods, the finite element approximation functions are defined on a background mesh that does not conform to the boundaries of the domain or internal material interfaces.  As immersed finite element methods employ an unfitted background mesh, they enable the simulation of physical systems that exhibit a change of domain topology, such as the flow of blood past the heart valves between the four main chambers of the human heart \cite{kamensky2017projection,kamensky2017immersogeometric,kamensky2015immersogeometric}.  In order for classical finite element methods to be applied to such problems, the domain must be frequently remeshed, and each time the domain is remeshed, the solution field must be remapped between the old and new body-fitted meshes.  For both isogeometric and immersed finite element methods, strong enforcement of boundary and interface conditions is quite difficult due to the non-interpolatory nature of the primal field approximation space along the boundary of the domain. Thus, the weak Nitsche's approach is especially attractive in this context.  As an example, displacement and rotation boundary conditions are enforced weakly using Nitsche's method in the isogeometric Kirchhoff-Love shell analysis result displayed in Fig.~\ref{fig:shell_ex}, and continuity of velocity across the fluid-structure interface is weakly enforced using Nitsche's method in the immersed heart valve analysis result displayed in Fig.~\ref{fig:valve_ex}.

Unfortunately, the construction of Nitsche's method is problem-dependent and can be particularly tedious to construct for high-order PDEs.  To facilitate the process, we presented a general framework in \cite{Benzaken2020} for constructing Nitsche's method for weak enforcement of boundary conditions for variational problems. Given a generalized Green's identity, and suitable generalized trace and Cauchy-Schwarz inequalities, we showed how to derive Nitsche's method and established conditions under which the resulting method is both stable and convergent. In that original work, we were particularly interested in linearized Kirchhoff-Love shell discretizations. However, the framework is broadly applicable to any variational problem, although its application to new systems may still seem daunting. The goal of this paper is to clarify the process through a sequence of examples, starting with the familiar Poisson equation and extending it to vector-valued and higher-order PDEs.

It should be noted that Nitsche's method is closely related to discontinuous Galerkin methods of interior penalty type \cite{riviere2008discontinuous}.  In fact, for the Poisson equation, utilization of the symmetric interior penalty discontinuous Galerkin method for enforcing boundary conditions is equivalent to Nitsche's method \cite{arnold1982interior}.  It should come as no surprise, then, that the same steps that are required to construct Nitsche's method for a particular set of PDEs are also necessary to build symmetric interior penalty discontinuous Galerkin methods for the same set of PDEs.  However, we elect not to demonstrate this in the current tutorial to keep the tutorial focused on Nitsche's method.

The same steps that are required to construct Nitsche's method can also be used to build so-called non-symmetric Nitsche methods \cite{schillinger2016non}.  Non-symmetric Nitsche methods harbor certain advantages over standard Nitsche methods.  In particular, non-symmetric Nitsche methods do not require a penalty term to enforce coercivity.  They do, however, suffer from certain drawbacks, such as non-optimal convergence in lower-order Sobolev norms \cite{burman2012penalty}.  Again, to keep the tutorial focused and brief, we do not discuss non-symmetric Nitsche methods further here.

The remainder of this paper proceeds as follows. In Section~\ref{sec:Nitsche_construct}, we derive Nitsche's method for Poisson's scalar-valued equation, utilizing a constructive approach in efforts to build intution behind the formulation; we deliberately sacrifice mathematical rigor to establish intuition. Section~\ref{sec:Nitsche} presents the abstract framework for constructing Nitsche's method for variational problems originally introduced in \cite{Benzaken2020}. In Section~\ref{sec:Poisson}, we elucidate the concepts discussed in Section~\ref{sec:Nitsche} by returning to the basic Poisson problem and extending it to the vector-valued setting. Next, in Section~\ref{sec:biharmonic}, the framework is applied to the biharmonic equation. Since this is a fourth-order equation, multiple boundary conditions must be handled on each part of the boundary. Finally, Section~\ref{sec:KL_plate} extends the results from the biharmonic equation to the linearized Kirchhoff-Love plate equation. This requires a transformation of the physical boundary conditions, as well as the requirement of a slightly smoother solution, to fit within the framework for Nitsche's method. Section~\ref{sec:conclusion} provides some concluding remarks.


\section{Building intuition for Nitsche's method}
\label{sec:Nitsche_construct}

We begin our discussion of Nitsche's method by considering the Poisson equation. This includes estimation of the various constants and proving coercivity of the resulting bilinear form. As mentioned above, we sacrifice some mathematical rigor in this section for didactic purposes. The full rigor is included in the abstract framework that is reviewed in the next section. Throughout this section, the bold-faced phrases are the necessary ingredients for the abstract framework discussed in Section~\ref{sec:Nitsche}.

\subsection{The Scalar Poisson Problem}

Let $\physDomain \subset \R^\sDim$ be a bounded domain (i.e., a bounded, open set) with Lipschitz-continuous boundary $\physBoundary = \partial \physDomain$, where $\sDim \in \N$ is the spatial dimension. Furthermore, let $\mathcal{V}^{\SM} \equiv H^1 (\physDomain)$, $\textup{f} \in L^2(\physDomain)$, $\textup{g} \in H^{1/2}(\physBoundary)$, and define $\mathcal{V}^{\SM}_{\textup{g}} \equiv \left\{ v \in \mathcal{V}^{\SM}: v|_{\physBoundary} = \textup{g} \right\}$, where $v|_{\physBoundary}$ denotes the trace of $v \in \mathcal{V}^{\SM}(\physDomain)$ on $\physBoundary$. We are interested in the following strong form of the Poisson boundary value problem:

$$
(S^{\SM}) \left\{ \hspace{5pt}
\parbox{4.35in}{
\noindent \textup{Find $u: \overline{\physDomain} \rightarrow \R$ such that}
\begin{equation*}
\begin{aligned}\begin{array}{rll}
- \Delta u &= \textup{f} \hspace{10pt} &\textup{in} \ \physDomain\\
u &= \textup{g} \hspace{10pt} &\textup{on} \ \physBoundary.\\
\end{array}
\end{aligned}
\label{eqn:SM_Strong}
\end{equation*}
}
\right.
$$

We use superscript $\SM$ to denote the \emph{heuristic} Poisson problem discussed in this section to distinguish it from the other problems considered later in this paper. The first step towards a Nitsche formulation is to construct the corresponding variational problem, to which we then employ a Bubnov-Galerkin method to arrive at a discretization.

\subsection{The Variational Formulation}

To construct the variational problem, we first identify a Sobolev space in which a weak solution to this PDE must live. In this case, the solution $u \in \mathcal{V}^{\SM}_{\textup{g}}$. We designate the homogeneous counterpart to this space as the test space, and multiply both sides by a test function, $v \in \mathcal{V}^{\SM}_0$. Integration by parts on the left hand side using \textbf{Green's first identity} then gives
\begin{equation}
 - \int_\physDomain v \Delta u \ d \physDomain = \int_\physDomain \nabla u \cdot \nabla v \ d \physDomain - \int_\physBoundary \left( \nabla u \cdot {\bf n} \right) v \ d \physBoundary,
  \label{eqn:IBP_model}
\end{equation}
where ${\bf n}$ is the outward-facing unit normal to $\Omega$.

\begin{remark}
  In arriving at Green's first identity, we used the fact that $u$ is the solution to $(S^{\SM})$. That is, $u$ has (at least) a continuous Laplacian and therefore all integrals appearing in \eqref{eqn:IBP_model} are well-defined. However, this Green's identity is valid for a larger class of functions. In particular, let
\begin{equation*}
		\tilde{\mathcal{V}}^{\SM} \equiv \left\{ w \in H^1(\physDomain) : \Delta w \in L^2(\physDomain) \textup{ and } \nabla w \cdot {\bf n} \big|_{\physBoundary} \in L^2(\physBoundary) \right\},
\end{equation*}
then \eqref{eqn:IBP_model} holds with $u$ replaced by $w \in \tilde{\mathcal{V}}^{\SM}$ and $v \in \mathcal{V}^{\SM}$. As discussed later in Section~\ref{sec:Nitsche}, we will need spaces such as the one presented in this remark to arrive at a Nitsche formulation through our abstract framework.

\end{remark}

Importantly, the boundary integral on the right-hand side of \eqref{eqn:IBP_model} vanishes because $v \in \mathcal{V}^{\SM}_0$, i.e., $v \equiv 0$ on $\physBoundary$. The remaining term on the right hand side gives rise to the following bilinear form for all $w \in \mathcal{V}^{\SM}_{\textup{g}}$ and $v \in \mathcal{V}^{\SM}_0$:

\begin{equation}
a^{\SM}(w,v) \equiv \int_\physDomain \nabla w \cdot \nabla v \ d\physDomain.
\label{eqn:aM_cons}
\end{equation}
Similarly, integration of the forcing function against a test function yields the linear functional given by \begin{equation}
\left\langle \linFunctional^{\SM}, v \right\rangle \equiv \int_\physDomain \textup{f} v \ d\physDomain
\label{eqn:fM_cons}
\end{equation}
for all $v \in \mathcal{V}^{\SM}_0$. The definitions of \eqref{eqn:aM_cons} and \eqref{eqn:fM_cons} naturally lead to the following variational formulation:

$$
(V^{\SM}) \left\{ \hspace{5pt}
\parbox{4.35in}{
\noindent Find $u \in \mathcal{V}^{\SM}_{\textup{g}}$ such that
\begin{eqnarray*}
a^{\SM}(u, v) = \left\langle \linFunctional^{\SM}, v \right\rangle
\end{eqnarray*}
for every $v \in \mathcal{V}^{\SM}_0$.
}
\right.
$$

\begin{remark}
The variational problem $(V^{\SM})$ corresponds to the first order optimality conditions of the following minimization problem:
$$
(M^{\SM}) \left\{ \hspace{5pt}
\parbox{4.35in}{
\noindent Given $\textup{f}$ and $\textup{g}$, find $u \in \mathcal{V}^{\SM}_{\textup{g}}$ that minimizes the total energy
\begin{eqnarray*}
E^{\SM}_{\text{total}}(u) = \frac{1}{2} \int_\physDomain \nabla u \cdot \nabla u \ d\physDomain - \int_\physDomain \textup{f} u \ d\physDomain.
\end{eqnarray*}
}
\right.
$$
Note that $(V^{\SM})$ is obtained by setting the first variation of the energy appearing in $(M^{\SM})$ to zero. In fact, this is where the term ``variational form'' comes from. Moreover, the abstract framework for constructing Nitsche's method presented in Section~\ref{sec:Nitsche} is grounded in this type of minimization problem, as we will discuss further in that section.
\end{remark}

At this point, we are able to select discrete subspaces of the trial space, $\mathcal{V}^{\SM}_{\textup{g}}$, and the test space, $\mathcal{V}^{\SM}_0$, to arrive at a Bubnov-Galerkin discretization of $(V^\SM)$. To accurately pose the associated discrete variational problem, we need to define a mesh over which the problem is discretized.

\subsection{Bubnov-Galerkin Method}

Assume that $\Omega$ can be represented in terms of a mesh $\mathcal{K}$ of non-overlapping (mapped) polygons, i.e., elements, so that $\Omega = \textup{int}\left( \overline{\cup_{K \in \mathcal{K}} K}\right)$. Assume further that the approximation space $\mathcal{V}^{\SM}_h \subset \mathcal{V}^\SM$ consists of (at least) $C^0$-continuous piecewise polynomial or rational approximations over the mesh $\mathcal{K}$. For each element $K \in \mathcal{K}$, we associate an element size $h_{K} = \textup{diam}(K)$ and we associate with the entire mesh $\mathcal{K}$ the mesh size $h = \max_{K \in \mathcal{K}} h_K$.

With a mesh in hand, we can define the discrete trial space as the subspace of $\mathcal{V}^{\SM}_h$ that satisfies the Dirichlet boundary conditions on the edge mesh, i.e., $\mathcal{V}^{\SM}_{\textup{g},h} = \mathcal{V}^{\SM}_{\textup{g}} \cap \mathcal{V}^{\SM}_h$, and the discrete test space as the homogeneous counterpart. The discrete variational problem is then given by
$$
(V^{\SM}_h) \left\{ \hspace{5pt}
\parbox{4.35in}{
\noindent Find $u_h \in \mathcal{V}^{\SM}_{\textup{g},h}$ such that
\begin{eqnarray*}
a^{\SM}(u_h, v_h) = \left\langle \linFunctional^{\SM}, v_h \right\rangle
\end{eqnarray*}
for every $v_h \in \mathcal{V}^{\SM}_{0,h}$.
}
\right.
$$

\begin{remark}
The variational problem $(V_h^{\SM})$ corresponds to the first order optimality conditions of the following minimization problem:
$$
(M_h^{\SM}) \left\{ \hspace{5pt}
\parbox{4.35in}{
\noindent Given $\textup{f}$ and $\textup{g}$, find $u_h \in \mathcal{V}^{\SM}_{\textup{g},h}$ that minimizes the total energy
\begin{eqnarray*}
E^{\SM}_{h}(u_h) = \frac{1}{2} \int_\physDomain \nabla u_h \cdot \nabla u_h \ d\physDomain - \int_\physDomain \textup{f} u_h \ d\physDomain.
\end{eqnarray*}
}
\right.
$$
\end{remark}

There are three important properties associated with Galerkin's method that are central to this paper. In the context of the model problem considered in this section, these three properties are:
\begin{enumerate}
\vspace{.5em}
\item{
\textbf{Symmetry}: For $w_h, v_h \in \mathcal{V}_{0,h}$,
\begin{equation*}
	a^{\SM}(w_h,v_h) = a^{\SM}(v_h,w_h).
\end{equation*}
A symmetric bilinear form yields a symmetric linear system after discretization and is crucial in proving optimal asymptotic discretization error convergence rates in low-ordered Sobolev norms \cite[Thm. 3.7]{Strang1973}.
}

\vspace{.5em}
\item{
\textbf{Coercivity}: For $w_h \in \mathcal{V}_{\text{g},h}$
\begin{equation*}
	a^{\SM}(w_h,w_h) > 0,
\end{equation*}
for $w_h \not\equiv 0$. The coercivity condition is used for establishing uniqueness of the solution. This is accomplished via the Lax-Milgram theorem \cite[\S 6.2]{EvansPDEs}, which ensures the existence of a unique solution, provided that the bilinear form is both bounded and coercive.
}

\vspace{.5em}
\item{
\textbf{Consistency}: For the solution, $u$, to $(V^\SM)$
\begin{equation*}
	a^{\SM}(u, v_h) - \left\langle \linFunctional^{\SM}, v_h \right\rangle = 0.
\end{equation*}
for all $v_h \in \mathcal{V}^{\SM}_{0,h}$. Consistency is synonymous to \emph{Galerkin Orthogonality}, which appears frequently in literature. The two concepts are related via
\begin{equation*}
	\underbrace{ 0 = a^{\SM}(u - u_h, v_h) }_{\text{Galerkin Orthogonality}} = \underbrace{ a^{\SM}(u, v_h) - a^{\SM}(u_h, v_h) = 0 }_{\text{Consistency}}.
\end{equation*}
Therefore, a consistent discrete solution can be interpreted as one with error, $u - u_h$, that is $a^\SM$-orthogonal to the finite-dimensional subspace $\mathcal{V}^{\SM}_{0,h}$.
}

\end{enumerate}

Note that the discrete spaces $\mathcal{V}^{\SM}_{\textup{g},h}$ and $\mathcal{V}^{\SM}_{0,h}$ require strongly-enforced Dirichlet boundary conditions. This is straightforward for simple approximation spaces (e.g., piecewise linear finite elements), simple applications (e.g., Poisson's equation considered here), and simple constraints (e.g., displacement boundary conditions). However, for complex approximation spaces (e.g., B-splines and subdivision surfaces), complex applications (e.g., Kirchhoff-Love shells), and complex constraints (e.g., rotation boundary conditions for shells and plates), this becomes much more difficult to enforce. For these reasons, it is oftentimes more convenient to enforce Dirichlet boundary conditions weakly.

\subsection{Penalty Method}
\label{sec:penaltymethod}

Perhaps the simplest and most familiar approach to weakly enforcing boundary conditions is the penalty method. In this approach, we relinquish the strong boundary condition enforcement by expanding the test and trial spaces to include all of $\mathcal{V}^{\SM}_h$. However, \eqref{eqn:aM_cons} is not bounded and coercive on $\mathcal{V}_h^{\SM}$, as it was on $\mathcal{V}^{\SM}_{0,h}$ and therefore the necessary conditions of the Lax-Milgram theorem are not met. This is because the kernel of \eqref{eqn:aM_cons} over $\mathcal{V}^{\SM}_h$ is non-trivial in that it contains the set of all constant functions. We can restore coercivity by adding a residual-based {\color{Orchid} \textbf{penalty term}} to both the bilinear and linear forms. In particular, we modify the finite-dimensional counterparts of \eqref{eqn:aM_cons} and \eqref{eqn:fM_cons} to be
\begin{equation}
a^{\SM}_{\text{pen}}(w_h,v_h) \equiv \int_\physDomain \nabla w_h \cdot \nabla v_h \ d\physDomain {\color{Orchid} \ \underbrace{ + \frac{\cPen{}}{h} \int_\physBoundary w_h v_h \ d \physBoundary }_\text{Penalty Term} }
\label{eqn:aM_weak_pen}
\end{equation}
and
\begin{equation}
\left\langle \linFunctional^{\SM}_{\text{pen}}, v_h \right\rangle \equiv \int_\physDomain \textup{f} v_h \ d\physDomain {\color{Orchid} \ \underbrace{ + \frac{\cPen{}}{h} \int_\physBoundary \textup{g} v_h \ d \physBoundary }_\text{Penalty Term} }.
\label{eqn:fM_weak_pen}
\end{equation}
Here, we have used a purple color coding to highlight the penalty term and to differentiate it from the other terms that will be presented shortly. This theme of color coding the terms that are ultimately related to Nitsche's method is used throughout the remainder of the paper. The explicit $h$-dependence appearing in the penalty terms of \eqref{eqn:aM_weak_pen} and \eqref{eqn:fM_weak_pen} comes from a dimensionality argument made later in Section~\ref{sec:m_coer_and_penalty}.

Provided $\cPen{} > 0$, then \eqref{eqn:aM_weak_pen} restores coercivity to the bilinear form over the larger space $\mathcal{V}^{\SM}_h$ by introducing an additional ``energy'' corresponding to enforcing the boundary condition. The constant $\cPen{}$ is known as the \textbf{penalty parameter}. This parameter is typically chosen heuristically or experimentally by the importance of enforcing the Dirichlet boundary conditions. As $\cPen{} \rightarrow \infty$, the penalty method recovers strong satisfaction of the Dirichlet boundary conditions. However, a large penalty parameter relative to the other terms in the formulation leads to ill-conditioning of the resulting matrix system after discretization, and consequently hurts numerical performance. The discrete variational penalty problem is given by:
$$
(V_\text{pen}^{\SM}) \left\{ \hspace{5pt}
\parbox{4.35in}{
\noindent Find $u_{h} \in \mathcal{V}^{\SM}_h$ such that
\begin{eqnarray*}
a^{\SM}_{\text{pen}}(u_h, v_h) = \left\langle \linFunctional^{\SM}_{\text{pen}}, v_h \right\rangle
\end{eqnarray*}
for every $v_h \in \mathcal{V}^{\SM}_h$.
}
\right.
$$

\begin{remark}
Alternatively, the penalty method can be understood as a modification to a discrete energy minimization problem, $(M^{\SM}_h)$, by the supplemental energy associated with the accuracy of boundary condition enforcement. In particular,
$$
(M_{\textup{pen}}^{\SM}) \left\{ \hspace{5pt}
\parbox{4.35in}{
\noindent Given $\textup{f}$, $\textup{g}$, and $\cPen{} > 0$, find $u_h \in \mathcal{V}^{\SM}_h$ that minimizes the total energy
\begin{eqnarray*}
E^{\SM}_{\textup{pen}}(u_h) = \frac{1}{2} \int_\physDomain \nabla u_h \cdot \nabla u_h \ d\physDomain - \int_\physDomain \textup{f} u_h \ d\physDomain + \frac{\cPen{}}{2h} \int_\physBoundary ( u_h - g )^2 \ d \physBoundary.
\end{eqnarray*}
}
\right.
$$
Note, however, that in general the minimizer of $(M_{\textup{pen}}^{\SM})$ does not coincide with that of $(M^{\SM}_h)$.
\end{remark}

Although we have a coercive bilinear form, we have lost consistency with the original problem. To see this, suppose that $u$ is smooth enough so that \eqref{eqn:IBP_model} holds and observe that the bilinear form is inconsistent since
\begin{equation}
\begin{aligned}
a^{\SM}_{\text{pen}} & (u, v_h) - a^{\SM}_{\text{pen}} (u_h, v_h)\\
&= a^{\SM}_{\text{pen}} (u, v_h) - \left\langle \linFunctional^{\SM}_{\text{pen}}, v_h \right\rangle\\
&= \int_\physDomain \nabla u \cdot \nabla v_h \ d\physDomain + \frac{\cPen{}}{h} \int_\physBoundary u v_h \ d \physBoundary - \int_\physDomain \textup{f} v_h \ d\physDomain - \frac{\cPen{}}{h} \int_\physBoundary \textup{g} v_h \ d \physBoundary\\
&= \int_\physBoundary \left( \nabla u \cdot {\bf n} \right) v_h \ d \physBoundary + \frac{\cPen{}}{h} \int_\physBoundary  \left( u - \textup{g} \right) v_h \ d \physBoundary - \int_\physDomain \left[ \textup{f} + \Delta u \right] v_h \ d \physDomain\\
&= \int_\physBoundary \left( \nabla u \cdot {\bf n} \right)  v_h \ d \physBoundary\\
&\neq 0.
\end{aligned}
\label{eqn:pen_not_con}
\end{equation}
Here, the definitions \eqref{eqn:aM_weak_pen} and \eqref{eqn:fM_weak_pen} were used in the second equality and the residual integrals going from the third equality to the fourth vanish because $u_h = \text{g}$ a.e. on $\physBoundary$ and $\textup{f} = - \Delta u_h$ a.e. in $\physDomain$.

As a consequence of this inconsistency, one should expect arrested convergence rates after numerical discretization for reasonably-sized penalty parameters~\cite{FernandezMendez2004}. To rectify this issue, we seek to restore the consistency of the variational formulation, which is the crux of Nitsche's method.

\subsection{Establishing Consistency and Restoring Symmetry}
\label{sec:con_and_sym_M}

The remaining boundary integral in \eqref{eqn:pen_not_con} is responsible for the inconsistency of the penalty formulation. Therefore, supplementing the penalty bilinear form \eqref{eqn:aM_weak_pen} with this boundary integral, but of opposite sign, should restore the variational consistency. In fact, this is precisely the approach taken to arrive at the Nitsche formulation. By adding a boundary term that negates the effects of the boundary integral present in \eqref{eqn:IBP_model}, which we refer to as the {\color{Cerulean} \textbf{consistency term}}, we arrive at the following bilinear form defined for $w_h, v_h \in \mathcal{V}_h^\SM$:

\begin{equation}
a^{\SM}_{\text{con}}(w_h,v_h) \equiv \int_\physDomain \nabla w_h \cdot \nabla v_h \ d\physDomain {\color{Cerulean} \ \underbrace{ - \int_\physBoundary \left( \nabla w_h \cdot {\bf n} \right) v_h \ d \physBoundary }_\text{Consistency Term} } {\color{Orchid} \ \underbrace{ + \frac{\cPen{}}{h} \int_\physBoundary w_h v_h \ d \physBoundary }_\text{Penalty Term} }.
\label{eqn:aM_con}
\end{equation}

By the same steps as in \eqref{eqn:pen_not_con}, one can show that we have constructed a consistent bilinear form in \eqref{eqn:aM_con}. Later in Section~\ref{sec:m_coer_and_penalty}, we will show that this new bilinear form is also coercive over $\mathcal{V}_h^\SM$, provided that $\cPen{}$ is large enough. However, we have lost the symmetry that was afforded by the initial, Bubnov-Galerkin formulation \eqref{eqn:aM_cons} and the penalty formulation \eqref{eqn:aM_weak_pen}. A discretization utilizing \eqref{eqn:fM_weak_pen} and \eqref{eqn:aM_con} would in general observe the theoretically optimal convergence rates in the energy norm. However, due to the lack of symmetry, sub-optimal convergence rates in lower-ordered Sobolev norms may be observed. Furthermore, such a discretization would not yield a symmetric linear system and would therefore not enjoy the associated numerical benefits.

Inspired by the construction of the penalty terms, the symmetry of the bilinear form can be restored by adding a residual-based boundary integral that comprises the {\color{ForestGreen} \textbf{symmetry term}}. This term is constructed from the symmetric counterpart to the consistency term. When added to both the existing bilinear form and the linear form, we obtain

\begin{equation}
a^{\SM}_h (w,v) \equiv \int_\physDomain \nabla w \cdot \nabla v \ d\physDomain {\color{Cerulean} \ \underbrace{ - \int_\physBoundary \left( \nabla w \cdot {\bf n} \right) v \ d \physBoundary }_\text{Consistency Term} } {\color{ForestGreen} \ \underbrace{ - \int_\physBoundary \left( \nabla v \cdot {\bf n} \right) w \ d \physBoundary }_\text{Symmetry Term} } {\color{Orchid} \ \underbrace{ + \frac{\cPen{}}{h} \int_\physBoundary w v \ d \physBoundary }_\text{Penalty Term} }
\label{eqn:aM_Nit}
\end{equation}
and
\begin{equation}
\left\langle \linFunctional^{\SM}_h, v \right\rangle \equiv \int_\physDomain \textup{f} v \ d\physDomain {\color{ForestGreen} \ \underbrace{ - \int_\physBoundary \left( \nabla v \cdot {\bf n} \right) \textup{g} \ d \physBoundary }_\text{Symmetry Term} } {\color{Orchid} \ \underbrace{ + \frac{\cPen{}}{h} \int_\physBoundary \textup{g} v \ d \physBoundary }_\text{Penalty Term} }.
\label{eqn:fM_Nit}
\end{equation}
The consistency of \eqref{eqn:aM_Nit} is by construction and the symmetry thereof is apparent. The next subsection is dedicated to establishing the coercivity of $a^{\SM}_h (\cdot, \cdot)$.

\subsection{Ensuring Coercivity through Penalty Parameter Selection}
\label{sec:m_coer_and_penalty}

In Section~\ref{sec:con_and_sym_M}, we claimed that the additional penalty terms will restore coercivity of the bilinear form in \eqref{eqn:aM_Nit}. In this subsection, we discuss the process of selecting the associated penalty parameter, guided by a set of inequalities.

The first inequality is the familiar \textbf{Cauchy-Schwarz inequality} that is of the form
\begin{equation*}
  \int_\physBoundary \left| \left( \nabla w_h \cdot {\bf n} \right) w_h \right| d \physBoundary \le \| \nabla w_h \cdot {\bf n} \|_{0,\physBoundary} \| w_h \|_{0,\physBoundary}.
\end{equation*}
where
\begin{equation}
  \| f \|^2_{0,\physBoundary} \equiv \int_\physBoundary \left| f \right|^2 d \physBoundary,
\end{equation}
is the classical $L^2$-norm on $\physBoundary$.

In this instance, the Cauchy-Schwarz inequality is standard in that it bounds a norm of a product by a product of norms. However, it will later be employed to provide a bound on duality pairings that involve general boundary quantities.

The second type of inequality we need is a \textbf{discrete trace inequality} that is of the form
\begin{equation}
\| \nabla w_h \cdot {\bf n} \|^2_{0,\physBoundary} \le \frac{C_{\text{tr}}}{h} \| \nabla w_h \|^2_{0,\physDomain}.
\label{eqn:TI_model_cons}
\end{equation}
This relates normed boundary data to normed interior data through a constant and appropriate mesh scaling.

The value $C_{\text{tr}}$ is a trace inequality constant which is dependent on factors such as the geometry of the domain, $\physDomain$. However, it is often convenient to distill out the explicit dependence of the trace inequality constant on factors such as $h$, as we have done in \eqref{eqn:TI_model_cons}, as well as other problem-dependent coefficients. This can be quickly accomplished for the mesh dependency through the following observation:

\begin{equation}
  \overbrace{ \int_\physBoundary \underbrace{ \left| \nabla w_h \cdot {\bf n} \right|^2 \strut }_{ \strut \left( \mathcal{O}(h^{-1}) \right)^2} \underbrace{ \strut d \physBoundary }_{\strut \mathcal{O}(h^{\sDim-1})} }^{\| \nabla w_h \cdot {\bf n} \|^2_{0,\physBoundary}} \le \frac{C_{\text{tr}}}{h} \overbrace{ \int_\physDomain \underbrace{ \left| \nabla w_h \right|^2 \strut }_{\strut \left( \mathcal{O}(h^{-1}) \right)^2} \underbrace{ \strut d \physDomain }_{\strut \mathcal{O}(h^{\sDim})} }^{ \| \nabla w_h \|^2_{0,\physDomain} }.
\end{equation}

\noindent Therefore, $\| \nabla w_h \cdot {\bf n} \|^2_{0,\physBoundary} = \mathcal{O}(h^{\sDim - 3})$ and $\| \nabla w_h \|^2_{0,\physDomain} = \mathcal{O}(h^{\sDim - 2})$, so the additional factor of $h^{-1}$ must be included on the right-hand side to rectify this mesh scaling discrepancy. Further elaboration on trace inequalities, and a potential source for explicit values or means for estimation of such constants, can be found in \cite{Evans2013}. Note that this type of inequality can be generalized, as is done in later sections, to accommodate the boundary conditions that are encountered for more exotic PDEs and those of a higher order.

Given these inequalities, we now arrive at a condition to ensure the coercivity of $a_{h}^{\SM}(\cdot,\cdot)$. Let $\CSconstant > 0$ and observe that

\begin{equation*}
\begin{aligned}
a_{h}^{\SM}(w_h,w_h) &= \int_\physDomain \nabla w_h \cdot \nabla w_h \ d\physDomain - 2 \int_\physBoundary \left( \nabla w_h \cdot {\bf n} \right) w_h \ d \physBoundary + \frac{\cPen{}}{h} \int_\physBoundary w_h^2 \ d \physBoundary\\
&= \| \nabla w_h \|^2_{0,\physDomain} - 2 \int_\physBoundary \left( \nabla w_h \cdot {\bf n} \right) w_h \ d \physBoundary + \frac{\cPen{}}{h} \| w_h \|^2_{0,\physBoundary}\\
&\ge \| \nabla w_h \|^2_{0,\physDomain} - 2 \| \nabla w_h \cdot {\bf n} \|_{0,\physBoundary} \| w_h \|_{0,\physBoundary} + \frac{\cPen{}}{h} \| w_h \|^2_{0,\physBoundary}\\
&\ge \| \nabla w_h \|^2_{0,\physDomain} - 2\left( \frac{1}{2\CSconstant} \| \nabla w_h \cdot {\bf n} \|^2_{0,\physBoundary} + \frac{\CSconstant}{2} \| w_h \|^2_{0,\physBoundary} \right) + \frac{\cPen{}}{h} \| w_h \|^2_{0,\physBoundary}\\
&\ge \| \nabla w_h \|^2_{0,\physDomain} - \left( \frac{C_{\text{tr}}}{h \CSconstant} \| \nabla w_h \|^2_{0,\physDomain} + \CSconstant \| w_h \|^2_{0,\physBoundary} \right) + \frac{\cPen{}}{h} \| w_h \|^2_{0,\physBoundary}\\
&= \left(1 - \frac{C_{\text{tr}}}{h \CSconstant} \right)\| \nabla w_h \|^2_{0,\physDomain} + \left( \frac{\cPen{}}{h} - \CSconstant \right) \| w_h \|^2_{0,\physBoundary}.
\end{aligned}
\end{equation*}

\noindent Note that we have used a Cauchy-Schwarz inequality going from line 2 to line 3, Young's inequality going from line 3 to line 4, and the trace inequality going from line 4 to line 5. We did not specifically present Young's inequality before the coercivity proof because it is a standard inequality that is not problem-dependent. This will also be the case for the inverse inequalities that appear later in the context of higher order PDEs. With the above relationship, we can ensure the coercivity of the bilinear form, i.e. that $a_{h}^{\SM}(w_h,w_h) > 0$, by selecting $\CSconstant$ and $\cPen{}$ such that
\begin{equation*}
    \CSconstant > \frac{C_{\text{tr}}}{h} \hspace{15pt} \text{and} \hspace{15pt} \frac{\cPen{}}{h} > \CSconstant,
\end{equation*}
or simply,
\begin{equation*}
  \cPen{} > C_{\text{tr}} > 0.
\end{equation*}

In practice, the choice of $\cPen{} = 2 C_{\text{tr}}$ yields good numerical results, provided that $C_{\text{tr}}$ can be estimated. We will revisit the estimation of $C_{\text{tr}}$ in Section~\ref{sec:model_TI_CS}. Note that this proof also demonstrates that the penalty parameter, like the trace inequality, scales linearly with respect to the reciprocal of the mesh size. In general, penalty parameters associated with Nitsche formulations are mesh dependent.

\subsection{Nitsche's Method}

At this point we are finally in a position to state Nitsche's method for the Poisson problem:

\vspace{1em}
\begin{mybox}[\emph{Heuristic Nitsche's Method for the Poisson Problem}]
  \vspace{-3pt}
$$
\hspace{-5pt}
(N_h^{\SM}) \left\{ \hspace{5pt}
\parbox{4.35in}{
For $\cPen{} > 0$, find $u_h \in \mathcal{V}^{\SM}_h$ such that
\begin{align*}
  a_{h}^{\SM}(u_h,v_h) = \underbrace{ \int_\physDomain \textup{f} v_h \ d\physDomain }_{\left \langle \linFunctional^{\SM}, v_h \right \rangle} {\color{ForestGreen} \underbrace{ - \int_\physBoundary \left( \nabla v_h \cdot {\bf n} \right) \textup{g} \ d \physBoundary }_\text{Symmetry Term} } {\color{Orchid} \ \underbrace{ + \frac{\cPen{}}{h}  \int_\physBoundary  \textup{g} v_h \ d \physBoundary }_\text{Penalty Term}  }
 \end{align*}
for every $v_h \in \mathcal{V}^{\SM}_h$ where $a_{h}^{\SM}(\cdot,\cdot) \colon \mathcal{V}^{\SM}_h \times \mathcal{V}^{\SM}_h $ is the bilinear form defined via
\begin{equation*}
	\begin{aligned}
  a_{h}^{\SM}(w_h,v_h) &\equiv \underbrace{ \int_\physDomain \nabla w_h \cdot \nabla v_h \ d\physDomain }_{a^{\SM}(w_h,v_h)}  {\color{Cerulean} \underbrace{ - \int_\physBoundary \left( \nabla w_h \cdot {\bf n} \right) v_h \ d \physBoundary }_\text{Consistency Term} } \\
  &\phantom{\equiv} {\color{ForestGreen} \ \underbrace{ - \int_\physBoundary \left( \nabla v_h \cdot {\bf n} \right) w_h \ d \physBoundary }_\text{Symmetry Term} } {\color{Orchid} \ \underbrace{ + \frac{\cPen{}}{h} \int_\physBoundary  w_h v_h \ d \physBoundary }_\text{Penalty Term} } ,
  	\end{aligned}
\end{equation*}
for $w_h,v_h \in \mathcal{V}^{\SM}_h$.
}
\right.
\vspace{3pt}
$$
\end{mybox}
\vspace{1em}

We refer to the above formulation as the ``heuristic'' formulation because it is missing some subtle points that have been omitted for clarity. We will revisit this problem in Section~\ref{sec:Poisson} where we will provide these details.

\begin{remark}
Alternatively, Nitsche's method can be understood as a modification to the discrete penalty energy minimization problem, $(M_\textup{pen}^{\SM})$, by the supplemental energy associated with the variationally-consistent boundary condition enforcement. In particular,
$$
(M_{\textup{nit}}^{\SM}) \left\{ \hspace{5pt}
\parbox{4.35in}{
\noindent Given $\textup{f}$, $\textup{g}$, and $\cPen{} > 0$, find $u_h \in \mathcal{V}^{\SM}_h$ that minimizes the total energy
\begin{equation*}
	\begin{aligned}
E^{\SM}_{\textup{nit}}(u_h) &= \frac{1}{2} \int_\physDomain \nabla u_h \cdot \nabla u_h \ d\physDomain - \int_\physDomain \textup{f} u_h \ d\physDomain - \int_\physBoundary \left( \nabla u_h \cdot {\bf n} \right) ( u_h - g )  \ d \physBoundary \\
&\phantom{=} + \frac{\cPen{}}{2h} \int_\physBoundary ( u_h - g )^2 \ d \physBoundary.
	\end{aligned}
\end{equation*}
}
\right.
$$
\end{remark}


\section{Abstract framework for Nitsche's Method}
\label{sec:Nitsche}

Having derived a Nitsche formulation for a simple model problem with a heuristic approach, we revisit the methods used in our construction with more formalism in this section. In particular, we recall the result of the Nitsche framework derived in \cite{Benzaken2020}. For brevity, we omit the associated proofs and in-depth discussions, including the connection to stabilized Lagrange multipliers formulations; the interested reader is referred to the original paper for more details.

Let $\mathcal{V}$ and $\mathcal{Q}$ be two Hilbert spaces with respective inner products $(\cdot,\cdot)_\mathcal{V}$ and $(\cdot,\cdot)_\mathcal{Q}$ and induced norms $\| \cdot \|_\mathcal{V} = (\cdot,\cdot)^{1/2}_\mathcal{V}$ and $\| \cdot \|_\mathcal{Q} = (\cdot,\cdot)^{1/2}_\mathcal{Q}$. We also use the notation $| \cdot |$ to refer to the absolute value for scalar quantities and the Euclidean norm for vector quantities. Let $\mathcal{V}^*$ and $\mathcal{Q}^*$ be the respective dual spaces of $\mathcal{V}$ and $\mathcal{Q}$, let ${}_{\mathcal{V}^*}\langle \cdot, \cdot \rangle_{\mathcal{V}}$ be the duality pairing between $\mathcal{V}$ and its dual, and let ${}_{\mathcal{Q}^*}\langle \cdot, \cdot \rangle_{\mathcal{Q}}$ denote the duality pairing between $\mathcal{Q}$ and its dual. Furthermore, let $\mathcal{T}: \mathcal{V} \rightarrow \mathcal{Q}$ be a bounded, surjective linear map referred to as the trace operator. Then given $g \in \mathcal{Q}$, define
\begin{equation*}
\mathcal{V}_g := \left\{ v \in \mathcal{V}: \mathcal{T}v = g \right\}.
\end{equation*}
Finally, let $a : \mathcal{V} \times \mathcal{V} \rightarrow \mathbb{R}$ be a bounded, symmetric, positive semi-definite bilinear form satisfying the following coercivity condition on the kernel of $\mathcal{T}$:
\begin{equation*}
a(v,v) \geq C \| v \|^2_\mathcal{V} \hspace{10pt} \forall v \in \mathcal{V}_0
\label{eqn:coerA}
\end{equation*}
for some constant $C \in \mathbb{R}_+$.

We are interested in the following minimization problem:
$$
(M) \left\{ \hspace{5pt}
\parbox{4.35in}{
\noindent Given $f \in \mathcal{V}^*$ and $g \in \mathcal{Q}$, find $u \in \mathcal{V}_g$ that minimizes the total energy
\begin{eqnarray*}
E_{\textup{total}}(u) = E_{\textup{int}}(u) + E_{\textup{ext}}(u)
\end{eqnarray*}
where the \textbf{\textit{internal energy}} is defined by
\begin{eqnarray*}
E_{\textup{int}}(u) \equiv \frac{1}{2} a(u,u)
\end{eqnarray*}
and the \textbf{\textit{external energy}} is defined by
\begin{eqnarray*}
E_{\textup{ext}}(u) \equiv -{}_{\mathcal{V}^*}\langle f, u \rangle_{\mathcal{V}}.
\end{eqnarray*}
}
\right.
$$
Note that we define the external energy $E_{\textup{ext}}(u)$ to be the negative of the duality pairing, not the other way around. In general, one can construct a valid functional representing the external energy that does not immediately appear to be a duality pairing, as we will see in Section~\ref{sec:KL_plate}. However, in order for this functional to be amenable to our framework, it must be mathematically manipulated to take this appropriate form of a duality pairing. This notion of external work arising from a duality pairing can be interpreted physically: for every admissible displacement field, there is a corresponding, or dual, set of forces that induces the field. The negative sign signifies that these forces are doing work \emph{on} the system.

The principle of least action states that the variation of the total energy vanishes at the minimum. Consequently, Problem $(M)$ is equivalent to the following variational problem:
$$
(V) \left\{ \hspace{5pt}
\parbox{4.35in}{
\noindent Given $f \in \mathcal{V}^*$ and $g \in \mathcal{Q}$, find $u \in \mathcal{V}_g$ such that
\begin{eqnarray*}
a(u, \delta u) = {}_{\mathcal{V}^*}\langle f, \delta u \rangle_{\mathcal{V}}
\label{eq:virtual_work}
\end{eqnarray*}
for every $\delta u \in \mathcal{V}_0$.
}
\right.
$$
The Lax-Milgram theorem guarantees that Problem $(V)$ has a unique solution $u \in \mathcal{V}$ that depends continuously on the input data $f \in \mathcal{V}^*$ and $g \in \mathcal{Q}$ \cite{EvansPDEs}.

Let $\mathcal{V}_h \subset \mathcal{V}$ be a finite-dimensional approximation space and $\mathcal{V}_{g,h} = \mathcal{V}_h \cap \mathcal{V}_g$ for every $g \in \mathcal{Q}$.  The Bubnov-Galerkin approximation of Problem $(V)$ then reads as follows:
$$
(V_h) \left\{ \hspace{5pt}
\parbox{4.35in}{
\noindent Given $f \in \mathcal{V}^*$ and $g \in \mathcal{Q}$, find $u_h \in \mathcal{V}_{g,h}$ such that
\begin{eqnarray*}
a(u_h, \delta u_h) = {}_{\mathcal{V}^*}\langle f, \delta u_h \rangle_{\mathcal{V}}
\end{eqnarray*}
for every $\delta u_h \in \mathcal{V}_{0,h}$.
}
\right.
$$
The Lax-Milgram theorem also guarantees that Problem $(V_h)$ has a unique solution $u_h \in \mathcal{V}_h$ that depends continuously on the input data $f \in \mathcal{V}^*$ and $g \in \mathcal{Q}$, and it is easily shown that the solution to Problem $(V_h)$ best approximates the solution to Problem $(V)$ with respect to the norm induced by the bilinear form $a(\cdot,\cdot)$.  The difficulty associated with Problem $(V_h)$ is the need for strong enforcement of the condition $\mathcal{T}u_h = g$. Similar to our derivation in the previous section, we address this by resorting to Nitsche's method.

Before proceeding to the presentation of Nitsche's method for Problem $(V)$, we note that the following two assumptions must be satisfied to arrive at a formulation that is both provably stable and convergent. The first assumption pertains to the existence of a generalized Green's identity and a smoothness condition on the problem solution.

\begin{assumption}
The following two conditions hold:
  \begin{enumerate}

  \item[$(a)$]{
  There exists a dense subspace $\tilde{\mathcal{V}} \subset \mathcal{V}$ and linear maps $\mathcal{L} : \tilde{\mathcal{V}} \rightarrow \mathcal{V}^*$ and $\mathcal{B}: \tilde{\mathcal{V}} \rightarrow \mathcal{Q}^*$ such that the following \textit{\textbf{generalized Green's identity}} holds:
  \begin{equation}
      a(w,v) = {}_{\mathcal{V}^*}\langle \mathcal{L}w, v \rangle_{\mathcal{V}} + {}_{\mathcal{Q}^*}\langle \mathcal{B}w, \mathcal{T}v \rangle_{\mathcal{Q}}
      \label{eq:int_by_parts}
    \end{equation}
    for all $w \in \tilde{\mathcal{V}}$ and $v \in \mathcal{V}$.
    }

    \item[$(b)$]{
    The solution $u$ of Problem $(M)$ satisfies $\mathcal{L} u = f$ whenever $f \in \mathcal{V}^*$ and $g \in \mathcal{Q}$ are such that $u \in \tilde{\mathcal{V}}$.
    }

    \end{enumerate}
    \label{assumption1}
\end{assumption}

\noindent The existence of a generalized Green's identity allows us to relate quantities in the interior of $\physDomain$ to the boundary $\physBoundary$ through classic theoretical results from vector calculus. More specifically, it elucidates the variationally-consistent boundary terms that connect the strong and weak formulations.

\begin{remark}
  For most commonly encountered operators, $(b)$ follows immediately from $(a)$ in Assumption~\ref{assumption1}. However, it is possible for general operators that $(b)$ does not immediately follow from $(a)$ although the authors have not been able to construct a simple example. To be completely rigorous, one must therefore show that both parts of the assumption are satisfied for the particular problems of interest. For the problems considered in this tutorial, it is indeed the case that $(a)$ implies $(b)$, so for simplicity, we only prove $(a)$ and refer the reader to \cite{Benzaken2020} for the complete details of how $(b)$ is shown for a specific problem, namely a Kirchhoff-Love shell.
\label{rem:skip_Lu_is_f}
\end{remark}

\begin{remark}
  In practice, the generalized Green's identity in \eqref{eq:int_by_parts} is obtained by applying integration by parts to the original variational formulation. The resulting map $\mathcal{L}$ encodes the differential-algebraic operators associated with the governing system of PDEs in their strong form as well as those associated with the natural boundary conditions. Similarly, the map $\mathcal{B}$ encodes the energetically conjugate natural boundary conditions that result from the application of integration by parts. Note that in order for the generalized Green's identity to hold, the solution to Problem $(V)$ must be sufficiently smooth. This is the reason we introduced an additional subspace $\tilde{\mathcal{V}} \subset \mathcal{V}$ for which \eqref{eq:int_by_parts} holds.
  \label{remark:V_tilde}
\end{remark}

\begin{remark}
  To distinguish the Green's identity given in \eqref{eq:int_by_parts} from Green's first, second, and third identities, we have used the clarifier ``generalized''.  For the scalar Poisson problem considered in Section~\ref{sec:Nitsche_construct}, the generalized Green's identity given by \eqref{eq:int_by_parts} coincides with Green's first identity used in \eqref{eqn:IBP_model}.
\end{remark}

The second assumption pertains to the existence of generalized trace and Cauchy-Schwarz inequalities.

\begin{assumption}
  \label{assumption2}
  There exists densely defined, positive, self-adjoint linear maps $\eta: \textup{dom}(\eta) \subseteq \mathcal{Q}^* \rightarrow \mathcal{Q}$ and $\epsilon: \textup{dom}(\epsilon) \subseteq \mathcal{Q}^* \rightarrow \mathcal{Q}$ such that $\epsilon$ is surjective (and thus invertible) and the following properties are satisfied:
\begin{enumerate}
\item The domain of definition of the operator $\mathcal{B} \colon \tilde{\mathcal{V}} \rightarrow \mathcal{Q}^*$ can be extended to the enlarged space $\tilde{\mathcal{V}} + \mathcal{V}_h$ and the space $\left\{ \mathcal{B}v: v \in \tilde{\mathcal{V}} + \mathcal{V}_h \right\}$ is a subset of $\textup{dom}(\eta)$.
\item The \textit{\textbf{generalized trace inequality}}:
\begin{equation*}
{}_{\mathcal{Q}^*}\langle \mathcal{B}v_h, \eta \mathcal{B}v_h \rangle_{\mathcal{Q}} \leq a(v_h,v_h)
\label{eqn:gen_trace}
\end{equation*}
holds for all $v_h \in \mathcal{V}_h$.
\item The \textit{\textbf{generalized Cauchy-Schwarz inequality}}:
\begin{equation*}
\left| {}_{\mathcal{Q}^*}\langle \mathcal{B} v, \mathcal{T} w \rangle_{\mathcal{Q}} \right| \leq \frac{1}{\CSconstant} {}_{\mathcal{Q}^*}\langle \mathcal{B}v, \eta \mathcal{B}v \rangle^{1/2}_{\mathcal{Q}} {}_{\mathcal{Q}^*}\langle \epsilon^{-1} \mathcal{T} w, \mathcal{T} w \rangle^{1/2}_{\mathcal{Q}}
\label{eqn:gen_CS}
\end{equation*}
holds for all $v, w \in \tilde{\mathcal{V}} + \mathcal{V}_h$, where $\CSconstant \in (1,\infty)$.
\end{enumerate}
\end{assumption}

The notion of domain extension presented in Assumption~\ref{assumption2} is further expounded upon in Section~\ref{sec:domain_enlargement}.

\begin{remark}
  The concept of vector space addition that appears in Assumption~\ref{assumption2}.1 may be unfamiliar. This operation is defined via
  \begin{equation}
    \tilde{\mathcal{V}} + \mathcal{V}_h := \left\{ \tilde{w} + {w}_h \colon \tilde{w} \in \tilde{\mathcal{V}} \text{ and } {w}_h \in \mathcal{V}_h \right\},
  \end{equation}
  that is, $\tilde{\mathcal{V}} + \mathcal{V}_h$ is the set of all possible combinations ${w} = \tilde{w} + {w}_h$ of elements in either subspace. This can easily be confused with the union of subspaces, which is defined as $\tilde{\mathcal{V}} \cup \mathcal{V}_h := \left\{ {w} \colon {w} \in \tilde{\mathcal{V}} \text{ or } {w} \in \mathcal{V}_h \right\}$. In particular, it is crucial to note that $\tilde{\mathcal{V}} + \mathcal{V}_h$ is a proper linear subspace of $\mathcal{V}$ while the same cannot be said in general for $\tilde{\mathcal{V}} \cup \mathcal{V}_h$. This is because if $\tilde{w} \in \tilde{\mathcal{V}}$ and ${w}_h \in \mathcal{V}_h$, we have that $\tilde{w} \in \tilde{\mathcal{V}} \cup \mathcal{V}_h$ and ${w}_h \in \tilde{\mathcal{V}} \cup \mathcal{V}_h$ but we do not necessarily have that $\tilde{w} + {w}_h \in \tilde{\mathcal{V}} \cup \mathcal{V}_h$. Finally, the vector space sum should not be confused with the direct sum of vector spaces, denoted by $\tilde{\mathcal{V}} \oplus \mathcal{V}_h$. The direct sum of spaces is a special case of the vector space sum wherein each element can be written \emph{uniquely} as a linear combination of one entity from either space. In other words, if $\tilde{\mathcal{V}} \cap \mathcal{V}_h = \{ 0 \}$, then $\tilde{\mathcal{V}} + \mathcal{V}_h \equiv \tilde{\mathcal{V}} \oplus \mathcal{V}_h$, however this is not the case in general.
\end{remark}

\begin{remark}
  Both $\eta$ and $\epsilon$ can be interpreted as Riesz operators. The former is associated with a duality pairing pertaining to the boundary operator $\mathcal{B}$ while the latter arises from a stabilized Lagrange multipliers method as discussed in \cite[\S 2]{Benzaken2020}. Note that the use of $\epsilon^{-1}$ rather than $\epsilon$ also comes from that discussion. Specifically, the Nitsche formulation is the analytical determination of the Lagrange Multipliers field that is solved through a static condensation, hence the inversion of $\epsilon$. Intuitively speaking, these operators are necessary to form a duality pairing rather than an inner product between two entities from the same space, i.e., $\| \mathcal{B} v_h \|^2_{\mathcal{Q}^*} = \left( \mathcal{B} v_h, \mathcal{B} v_h \right)_{\mathcal{Q}^*}$ vs. ${}_{\mathcal{Q}^*}\langle \mathcal{B}v_h, \eta \mathcal{B}v_h \rangle_{\mathcal{Q}}$ for $v_h \in \mathcal{V}_h$. For simple problems, these operators can take the form of a constant, cf. Remark~\ref{rem:eps_eta_VM}, while for more complex problems, these operators cannot be described so concisely, e.g., \eqref{eqn:eta_B} and \eqref{eqn:eps_B}.
\end{remark}

\noindent The generalized trace and Cauchy-Schwarz inequalities provide a mechanism for providing bounds for boundary quantities that are used in proving convergence of our Nitsche formulation. With the two assumptions above in hand, we are ready to present Nitsche's method for the abstract variational problem given by Problem $(V)$.

\vspace{1em}
\begin{mybox}[\emph{Nitsche's Method for an Abstract Variational Problem}]
  \vspace{-3pt}
$$
(N_h) \left\{ \hspace{5pt}
\parbox{4.35in}{
Given $f \in \mathcal{V}^*$ and $g \in \mathcal{Q}$, find $u_h \in \mathcal{V}_h$ such that
\begin{equation*}
  a_h(u_h,\delta u_h) = {}_{\mathcal{V}^*}\langle f, \delta u_h \rangle_{\mathcal{V}} \ {\color{ForestGreen} \underbrace{ - {}_{\mathcal{Q}^*}\langle \mathcal{B} \delta u_h, g \rangle_{\mathcal{Q}} }_\text{Symmetry Term} } \ {\color{Orchid} \underbrace{ + {}_{\mathcal{Q}^*}\langle \epsilon^{-1} \mathcal{T} \delta u_h, g \rangle_{\mathcal{Q}} }_\text{Penalty Term} }
  \label{eq:nitsche}
\end{equation*}
for every $\delta u_h \in \mathcal{V}_h$, where $a_h: \left( \tilde{\mathcal{V}} + \mathcal{V}_h \right) \times \left( \tilde{\mathcal{V}} + \mathcal{V}_h \right) \rightarrow \mathbb{R}$ is the bilinear form defined by
\begin{equation*}
  a_h(w, v) \equiv a(w, v) {\color{Cerulean} \ \underbrace{ - {}_{\mathcal{Q}^*}\langle \mathcal{B} w, \mathcal{T}v \rangle_{\mathcal{Q}} }_\text{Consistency Term} } \ {\color{ForestGreen} \underbrace{ - {}_{\mathcal{Q}^*}\langle \mathcal{B} v, \mathcal{T} w \rangle_{\mathcal{Q}} }_\text{Symmetry Term} } \ {\color{Orchid} \underbrace{ + {}_{\mathcal{Q}^*}\langle \epsilon^{-1} \mathcal{T} v, \mathcal{T} w \rangle_{\mathcal{Q}} }_\text{Penalty Term} }
\end{equation*}
for all $w, v \in \tilde{\mathcal{V}} + \mathcal{V}_h$.
}
\right.
\vspace{3pt}
$$
\end{mybox}
\vspace{1em}

Nitsche's method exhibits several important properties that give rise to its stability and convergence.  Namely, it is \textit{\textbf{consistent}}, its bilinear form $a_h(\cdot,\cdot)$ is \textit{\textbf{symmetric}}, and, provided the map $\epsilon$ is chosen appropriately, its bilinear form $a_h(\cdot,\cdot)$ is also \textit{\textbf{coercive}} on the discrete space $\mathcal{V}_h$. Before applying this framework to a set of examples, we present a fundamental result of Nitsche's method, namely, the well-posedness of the method as well as an error estimate.

\begin{theorem}[Well-Posedness and Error Estimate]
\label{theorem:error_estimate}
Suppose that Assumptions~\ref{assumption1} and~\ref{assumption2} hold.  Then there exists a unique discrete solution $u_h \in \mathcal{V}_h$ to Problem $(N_h)$.  Moreover, if the continuous solution $u \in \mathcal{V}$ to Problem $(M)$ satisfies $u \in \tilde{\mathcal{V}}$, then the discrete solution $u_h$ satisfies the error estimate
\begin{equation*}
\vvvertiii{u - u_h} \leq \left( 1+ \frac{2}{1-\frac{1}{\CSconstant}} \right) \min_{v_h \in \mathcal{V}_h} \vvvertiii{u - v_h},
\end{equation*}
where $\vvvertiii{\cdot}: \tilde{\mathcal{V}} + \mathcal{V}_h \rightarrow \mathbb{R}$ is the energy norm defined by
\begin{equation*}
  \vvvertiii{v}^2 \equiv a(v,v) + \hspace{-3pt}{\phantom{\big|}}_{\mathcal{Q}^*}\big\langle \mathcal{B} v, \eta \mathcal{B} v \big\rangle_{\mathcal{Q}} + 2 \hspace{-1pt}{\phantom{\big|}}_{\mathcal{Q}^*}\big\langle \left( \epsilon \right)^{-1} \mathcal{T} v, \mathcal{T} v \big\rangle_{\mathcal{Q}}
\end{equation*}
and $\gamma$ is given through Assumption~\ref{assumption2}.
\end{theorem}

Note that the above theorem applies to any formulation and problem setup for which Assumptions~\ref{assumption1} and~\ref{assumption2} hold.  Consequently, constructing Nitsche-based formulations for a new problem class should proceed according to the following steps:\\

\noindent \textbf{Step 1:} Construct an appropriate variational formulation, including specification of the Hilbert spaces $\mathcal{V}$ and $\mathcal{Q}$, the map $\mathcal{T}: \mathcal{V} \rightarrow \mathcal{Q}$, and the bilinear form $a : \mathcal{V} \times \mathcal{V} \rightarrow \mathbb{R}$.\\

\noindent \textbf{Step 2:} Establish the generalized Green's identity by determining the space $\tilde{\mathcal{V}}$ and the linear maps $\mathcal{L} : \tilde{\mathcal{V}} \rightarrow \mathcal{V}^*$ and $\mathcal{B} : \tilde{\mathcal{V}} \rightarrow \mathcal{Q}^*$ associated with Assumption~\ref{assumption1}. Note that the relevant operators will ultimately be defined over the extended domain $\tilde{\mathcal{V}} + \mathcal{V}_h$ for discretization.\\

\noindent \textbf{Step 3:} Establish the generalized Cauchy-Schwarz and trace inequalities by constructing suitable linear maps $\epsilon: \textup{dom}(\epsilon) \subseteq Q^* \rightarrow Q$ and $\eta: \textup{dom}(\eta) \subseteq Q^* \rightarrow Q$ such that Assumption~\ref{assumption2} is satisfied.\\

\noindent \textbf{Step 4:} Pose Nitsche's method according to Problem $(N_h)$.\\

In the following, we complete the above four steps to construct a Nitsche-based formulation for a vector-valued Poisson problem, a vector-valued biharmonic problem, and a linearized Kirchhoff-Love plate problem. Note that we do not need to conduct a full stability and convergence analysis, since we can readily employ the abstract framework presented here.  It should further be mentioned that symmetry can be employed to arrive at error estimates in norms other than the energy norm using the well-known Aubin-Nitsche trick \cite{ciarlet1991basic}.


\section{Nitsche's method for Poisson's Equation}
\label{sec:Poisson}

Using the abstract framework in the previous section, we proceed now to derive a Nitsche formulation for the vector-valued version of the Poisson model problem considered in Section~\ref{sec:Nitsche_construct}. The vector-valued setting presents a small increase in complexity that is easily handled by the abstract framework. In addition, we incorporate both Dirichlet and Neumann boundary conditions to demonstrate how the framework accommodates both types. However, the most significant difference compared to Section~\ref{sec:Nitsche_construct} is the increased mathematical rigor, which we highlight along the way. If necessary, we encourage the reader to refer back to Section~\ref{sec:Nitsche_construct} to recall the intuition that was built therein.


\subsection{The Variational Formulation}

Beginning with the first step in our recipe for constructing a Nitsche-based formulation, we first state the minimization problem and its corresponding variational formulation. This is often the most natural way to formulate a variational problem, even if it is not obvious in this particular example. Proceeding forward, we use superscript $\VM$ to refer to the vector-valued \emph{Poisson} problem considered here.

First, let $\physDomain \subset \R^\sDim$ be a bounded domain with Lipschitz-continuous boundary $\physBoundary = \partial \physDomain$, where $\sDim \in \N$ is the spatial dimension. Let ${\bf u}$, ${\bf v}$, and ${\bf w}$ denote vector-valued functions of dimension $\uDim \in \N$, where ${\bf u}$ is reserved to denote the solution, while ${\bf w}$ and ${\bf v}$ denote arbitrary trial and test functions, respectively, when they differ. We are interested in minimizing a functional with the following internal energy
\begin{eqnarray*}
E^{\VM}_{\textup{int}}({\bf w}) \equiv \frac{1}{2} \int_{\physDomain} \left( \nabla {\bf w} \right) : \left( \nabla {\bf w} \right) \ d\physDomain.
\end{eqnarray*}
Given this, we must select a space of admissible solutions that is smooth enough such that the functional is well-defined. We therefore restrict our attention to functions with at least \emph{one} integrable derivative. Thus, let $\mathcal{V}^{\VM} \equiv \left[ H^1(\physDomain) \right]^\uDim$ be the space of admissible solutions free of boundary conditions.

We partition the boundary into Dirichlet and Neumann parts, denoted $\physBoundary_D$ and $\physBoundary_N$, respectively. For a well-posed PDE, we require that $\physBoundary = \overline{\physBoundary_D \cup \physBoundary_N}$, $\physBoundary_D \cap \physBoundary_N = \emptyset$, and $\physBoundary_D \neq \emptyset$. With our choice of $\mathcal{V}^{\VM}$, the associated trace space is then $\mathcal{Q}^{\VM} \equiv \left[ H^{1/2}(\physBoundary_D) \right]^\uDim$

With the above boundary partitions in place, we are able to define the trace operator $\mathcal{T}^{\VM} \colon \mathcal{V}^{\VM} \rightarrow \mathcal{Q}^{\VM}$ via its action on an arbitrary ${\bf w} \in \mathcal{V}^{\VM}$ via $\mathcal{T}^{\VM} {\bf w} \equiv {\bf w} \big|_{\physBoundary_D}$, that is, ${\bf w}$ restricted to the Dirichlet boundary. The trace operator allows us to imbue our Hilbert space with boundary conditions, provided they are smooth enough, by the introduction of an additional constraint in the definition of the Hilbert space. In particular, given a prescribed displacement on the Dirichlet boundary, ${\bf g} \in \mathcal{Q}^{\VM}$, we denote the trial space of displacement fields satisfying the prescribed Dirichlet boundary conditions by
\begin{equation*}
\mathcal{V}^{\VM}_{{\bf g}} \equiv \left\{ {\bf w} \in \mathcal{V}^{\VM} \colon \mathcal{T}^{\VM} {\bf w} = {\bf g} \right\},
\end{equation*}
Similarly, $\mathcal{V}^{\VM}_{{\bf 0}}$ denotes the test space of displacement fields satisfying homogeneous Dirichlet boundary conditions, corresponding to ${\bf g} \equiv {\bf 0}$.

With only homogeneous forcing functions, the minimizer of $E^{\VM}_{\textup{int}}(\cdot)$ is dictated purely by the choice of the Dirichlet boundary condition ${\bf g}$. To accommodate the more general scenario with inhomogeneous forcing functions and/or Neumann-type boundary conditions, we additionally consider the functional representing external work that takes the form
\begin{eqnarray*}
E^{\VM}_{\textup{ext}}({\bf w}) \equiv - \int_{\physDomain} {\bf f} \cdot {\bf w} \ d \physDomain - \int_{\physBoundary_N} {\bf h} \cdot {\bf w} \ d \physBoundary.
\end{eqnarray*}
For this functional to be well-defined, we require that ${\bf f} \in \left[ L^2(\physDomain) \right]^\uDim$ and, for Neumann-type boundary conditions, that ${\bf h} \in \left[ L^2(\physBoundary_N) \right]^\uDim$.

We define the associated bilinear form $a^{\VM}(\cdot,\cdot): \mathcal{V}^{\VM} \times \mathcal{V}^{\VM} \rightarrow \mathbb{R}$ by taking the variation of the internal energy $E^{\VM}_{\textup{int}}({\bf w})$ with respect to the displacement field,
\begin{equation*}
a^{\VM}({\bf w},{\bf v}) \equiv \int_\physDomain \left( \nabla {\bf w} \right) : \left( \nabla {\bf v} \right) \ d\physDomain,
\label{eqn:aM}
\end{equation*}
for all ${\bf w}, {\bf v} \in \mathcal{V}^{\VM}$. Similarly, $\linFunctional^{\VM} \in \left( \mathcal{V}^{\VM} \right)^{\ast}$ is defined as the linear functional satisfying
\begin{equation*}
\left\langle \linFunctional^{\VM}, {\bf v} \right\rangle \equiv -E^{\VM}_{\textup{ext}}({\bf v})
\label{eqn:fM}
\end{equation*}
for all ${\bf v} \in \mathcal{V}^{\VM}$.

Given these definitions, the minimization problem of interest can be stated as follows:
$$
(M^{\VM}) \left\{ \hspace{5pt}
\parbox{4.35in}{
\noindent Find ${\bf u} \in \mathcal{V}^{\VM}_{{\bf g}}$ that minimizes the total energy
\begin{eqnarray*}
E^{\VM}_{\textup{total}}({\bf u}) = E^{\VM}_{\textup{int}}({\bf u}) + E^{\VM}_{\textup{ext}}({\bf u}).
\end{eqnarray*}
}
\right.
$$

The solution to Problem $(M^{\VM})$ is also the solution to the following variational problem:
$$
(V^{\VM}) \left\{ \hspace{5pt}
\parbox{4.35in}{
\noindent Find ${\bf u} \in \mathcal{V}^{\VM}_{{\bf g}}$ such that
\begin{eqnarray*}
a^{\VM}({\bf u}, \delta {\bf u}) = \left\langle \linFunctional^{\VM}, \delta {\bf u} \right\rangle \label{eqn:VM_Weak}
\end{eqnarray*}
for every $\delta {\bf u} \in \mathcal{V}^{\VM}_{\bf 0}$.
}
\right.
$$

\subsection{A Generalized Green's Identity}
\label{sec:greenpoisson}

The second step in our recipe for constructing a Nitsche-based formulation requires a generalized Green's identity that satisfies Assumption~\ref{assumption1} to be established. This step may in general not always be easy since it requires analytic integration by parts of the variational form. In this particular case, we begin the process by observing that, by the product rule,
\begin{equation}
  \nabla \cdot \left[ \left( \nabla {\bf w} \right) \cdot {\bf v} \right] = \left[ \nabla \cdot \left( \nabla {\bf w} \right) \right] \cdot {\bf v} + \left( \nabla {\bf w} \right) : \left( \nabla {\bf v} \right).
  \label{eqn:div_prod_rule_Poisson}
\end{equation}
Rearranging this expression, integrating both sides over the domain $\physDomain$, and applying the divergence theorem then yields
\begin{equation}
  \begin{aligned}
    a^{\VM}({\bf w},{\bf v}) &= \int_\physDomain \left( \nabla {\bf w} \right) : \left( \nabla {\bf v} \right) \ d \physDomain \\
    &= \int_\physDomain \nabla \cdot \left[ \left( \nabla {\bf w} \right) \cdot {\bf v} \right] \ d \physDomain - \int_\physDomain \left[ \nabla \cdot \left( \nabla {\bf w} \right) \right] \cdot {\bf v} \ d \physDomain\\
    &= \int_\physBoundary \left[ \left( \nabla {\bf w} \right) \cdot {\bf n} \right] \cdot {\bf v} \ d \physBoundary - \int_\physDomain {\bf v} \cdot \left( \Delta {\bf w} \right) \ d \physDomain,
  \end{aligned}
  \label{eqn:Greens_ID_SM_derive}
\end{equation}
where $\Delta {\bf w} \equiv \nabla \cdot \left(\nabla {\bf w}\right)$ is the Laplacian and ${\bf n}$ is a unit normal vector along $\physBoundary$ pointing outward from $\physDomain$. However, note that the generalized Green's identity arising from the last equality in \eqref{eqn:Greens_ID_SM_derive} is no longer valid for general ${\bf w} \in \mathcal{V}^{\VM}$ because it requires that (i) $\Delta {\bf w} \in \left[ L^2(\physDomain) \right]^\uDim$ and (ii) $\left( \nabla {\bf w} \right) \cdot {\bf n} \in \left[ L^2(\physBoundary) \right]^\uDim$. To satisfy these additional smoothness requirements, we let
\begin{equation*}
\tilde{\mathcal{V}}^{\VM} \equiv \left\{ {\bf v} \in \left[ H^1(\physDomain) \right]^\uDim : \Delta {\bf v} \in \left[L^2(\physDomain)\right]^\uDim \textup{ and } \left( \nabla {\bf v} \right) \cdot {\bf n} \big|_{\physBoundary} \in \left[L^2(\physBoundary)\right]^\uDim \right\}
\end{equation*}
and note that \eqref{eqn:Greens_ID_SM_derive} is now valid for general ${\bf w} \in \tilde{\mathcal{V}}^{\VM}$. Additionally, $\tilde{\mathcal{V}}^{\VM} \subset \mathcal{V}^{\VM}$ as required by Assumption~\ref{assumption1}.

\begin{remark}
  Our definition of $\tilde{\mathcal{V}}^{\VM}$ is the least-smooth space we can select, given the additional smoothness requirements elucidated by \eqref{eqn:Greens_ID_SM_derive}. However, it is not the only space that satisfies these newly-imposed requirements. For example, \eqref{eqn:Greens_ID_SM_derive} is also valid for general ${\bf w} \in \left[ H^2(\physDomain) \right]^\uDim$. The general properties in a classical sense of $\tilde{\mathcal{V}}^{\VM}$ are not immediately clear, but the Sobolev embedding theorem \cite{EvansPDEs} illustrates that $\left[ H^2(\physDomain) \right]^\uDim \subset \tilde{\mathcal{V}}^{\VM} \subset \left[ H^1(\physDomain) \right]^\uDim$, and for many problems the more restrictive choice of $\left[ H^2(\physDomain) \right]^\uDim$ will be sufficient. However, as an example, for domains with re-entrant corners, the solution ${\bf u} \in \mathcal{V}^{\VM}$ to Problem $(V^{\VM})$ does not generally lie in $\left[ H^2(\physDomain) \right]^\uDim$, while it does lie in $\tilde{\mathcal{V}}^{\VM}$ \cite{grisvard2011elliptic}.
\end{remark}


Although we have derived a generalized Green's identity, Assumption~\ref{assumption1} further stipulates that $\mathcal{L}^{\VM} {\bf u} = \linFunctional^{\VM}$ when ${\bf u} \in \tilde{\mathcal{V}}^{\VM}$. Typically at this juncture, this must be shown for complete rigor. However, as discussed in Remark~\ref{rem:skip_Lu_is_f}, we omit the details and simply assert that it is true.

We are then ready to state the following lemma regarding the generalized Green's identity for the vector-valued Poisson problem:

\begin{lemma}[Generalized Green's Identity for the Poisson Problem]
\label{lemma:greens_heat}
  For ${\bf w} \in \tilde{\mathcal{V}}^{\VM}$ and ${\bf v} \in \mathcal{V}^{\VM}$, the following generalized Green's identity holds:
  \begin{equation*}
    a^{\VM}({\bf w},{\bf v}) = \underbrace{ \int_{\physBoundary_N} \left[ \left( \nabla {\bf w} \right) \cdot {\bf n} \right]  \cdot {\bf v}\ d \physBoundary - \int_{\physDomain} (\Delta {\bf w}) \cdot {\bf v} \ d \physDomain }_{ \displaystyle \langle \mathcal{L}^{\VM} {\bf w}, {\bf v} \rangle } + \underbrace{ \int_{\physBoundary_D} \left[ \left( \nabla {\bf w} \right) \cdot {\bf n} \right] \cdot {\bf v} \ d\physBoundary }_{ \displaystyle \langle \mathcal{B}^{\VM} {\bf w}, \mathcal{T}^{\VM} {\bf v} \rangle }.
  \label{eqn:Greens_ID_M}
\end{equation*}
Moreover, the solution ${\bf u}$ of Problem $(V^{\VM})$ satisfies $\mathcal{L}^{\VM} {\bf u} = \linFunctional^{\VM}$ provided that the problem parameters ${\bf f}$, ${\bf g}$, and ${\bf h}$ are smooth enough such that ${\bf u} \in \tilde{\mathcal{V}}^{\VM}$.
  \begin{proof}
    The generalized Green's identity follows from the steps in \eqref{eqn:Greens_ID_SM_derive}. The result that $\mathcal{L}^{\VM} {\bf u} = \linFunctional^{\VM}$ follows by the same techniques used in \cite{Benzaken2020} as discussed in Remark~\ref{rem:skip_Lu_is_f}.
    
  \end{proof}
\end{lemma}

  \begin{remark}
    The corresponding strong form of Problem $(V^{\VM})$ is given by:
    $$
    (S^{\VM}) \left\{ \hspace{5pt}
    \parbox{4.35in}{
    \noindent \textup{Find ${\bf u}: \overline{\physDomain} \rightarrow \R$ such that:}
    \begin{equation*}
    \begin{aligned}\begin{array}{rll}
    - \Delta {\bf u} &= {\bf f} \hspace{10pt} &\textup{in} \ \physDomain\\
    {\bf u} &= {\bf g} \hspace{10pt} &\textup{on} \ \physBoundary_D\\
    \left( \nabla {\bf u} \right) \cdot {\bf n} &= {\bf h} \hspace{10pt} &\textup{on} \ \physBoundary_N\\
    \end{array}
    \end{aligned}
    \label{eqn:VM_Strong}
    \end{equation*}
    }
    \right.
    $$
    This result follows immediately from the relationship $\mathcal{L}^{\VM} {\bf u} = \linFunctional^{\VM}$ from Lemma~\ref{lemma:greens_heat} and the natural boundary conditions encoded by $\mathcal{B}^{\VM}$. Specifically, the linear functional $\linFunctional^{\VM}$ in Problem~$(V^{\VM})$ contains both the body forcing ${\bf f}$ and the natural boundary condition ${\bf h}$, while the variational space $\mathcal{V}_{\bf g}^\VM$, contains the essential boundary condition. On the other hand, the differential operator $\mathcal{L}^{\VM}$ encodes the energetically conjugate interior PDE and natural boundary condition to ${\bf f}$ and ${\bf h}$, respectively, while the boundary operator $\mathcal{B}^{\VM}$ encodes the energetically conjugate term associated with the essential boundary condition ${\bf g}$. Note that this is the same strong formulation, but in vector form, from which we started our derivations in Section~\ref{sec:Nitsche_construct}. This is readily seen by setting $\uDim = 1$.
  \end{remark}

\begin{remark}
A solution ${\bf u} \in \mathcal{V}^{\VM}$ to $(V^{\VM})$ is referred to as a weak solution to the minimization problem. When ${\bf u} \in \tilde{\mathcal{V}}^{\VM}$, the solution is smooth enough to satisfy the Euler-Lagrange equations associated with the minimization problem almost everywhere. For ${\bf u}$ to be a solution to $(S^{\VM})$, that is, a strong solution to the minimization problem, it must be even smoother. Specifically, by Sobolev embedding, ${\bf u}$ must be in a space endowed with the requisite number of continuous derivatives.

For example, in reference to the Poisson problem considered in this section, when $\textup{\bf f}$, $\textup{\bf g}$, and $\textup{\bf h}$ are such that ${\bf u} \in \left[ H^1 (\physDomain) \right]^\uDim$, then a weak solution to the problem exists, but it neither satisfies the Euler-Lagrange equations in a weak nor a strong sense, since the derivative operator is not defined pointwise. However, if $\textup{\bf f}$, $\textup{\bf g}$, and $\textup{\bf h}$ are such that ${\bf u} \in \left[ H^2 (\physDomain) \right]^\uDim$, then a weak solution exists and it satisfies the Euler-Lagrange equations almost everywhere, but it is not a strong solution, because $\Delta {\bf u}$ need not be well-defined everywhere. Finally, if $\textup{\bf f}$, $\textup{\bf g}$, and $\textup{\bf h}$ are such that ${\bf u} \in \left[ C^2 (\physDomain) \right]^\uDim$, then a weak solution exists and it satisfies the Euler-Lagrange equations both weakly and strongly, by Sobolev embedding.
  \label{rem:strong_weak_equiv}
\end{remark}

Following an identical procedure to that of Section~\ref{sec:Nitsche_construct}, we construct a mesh $\mathcal{K}$ associated with $\Omega$. We reiterate that our approximation space $\mathcal{V}^{\VM}_h$ consists of (at least) $C^0$-continuous piecewise polynomial or rational approximations over the mesh $\mathcal{K}$. The new ingredient pertaining to the mesh construction needed for this problem is that of the edge mesh. Specifically, each element $K \in \mathcal{K}$ has a set of edges, and we collect all of the edges of the mesh that lie along the boundary of the domain into an edge mesh $\mathcal{E}$. It is easily seen that each edge $E \in \mathcal{E}$ belongs to a unique element $K \in \mathcal{K}$. Thus, for each edge $E \in \mathcal{E}$, we associate an edge size $h_E = h_K$, where $K \in \mathcal{K}$ is the element for which $E$ is an edge. This is not the only size we can associate with the edge, but it is the simplest.  For anisotropic meshes, other prescriptions may be more appropriate (see, e.g., \cite{bazilevs2007weak}). We then further define the set of Dirichlet edges, that is, the set of edges belonging to the Dirichlet boundary, via
\begin{equation*}
  \mathcal{E}_D \equiv \left\{ E \in \mathcal{E} \colon E \subset \physBoundary_D \right\}.
\end{equation*}
We assume that $\physBoundary_D = \text{int}(\overline{\cup_{E\in\mathcal{E}_D} E})$ to ensure that each edge in $\mathcal{E}$ belongs to either the Dirichlet or the Neumann boundary, but not both.

\subsection{Domain Enlargement}
\label{sec:domain_enlargement}

In Section~\ref{sec:Nitsche_construct}, we showed consistency provided the solution was ``smooth enough'' and in Section~\ref{sec:Nitsche} we alluded to the fact that many of the operators involved in Nitsche formulations require a domain extension. In this subsection, we address these concepts formally.

The motivation for domain enlargement begins after the introduction of $\tilde{\mathcal{V}}^{\VM}$, where numerical practicality concerns may be raised due to the additional required smoothness. That is, a finite-dimensional subspace of $\tilde{\mathcal{V}}^{\VM}$ may require significantly more smoothness than that of a subspace of $\mathcal{V}^{\VM}$. We can relax this constraint through Assumption~\ref{assumption2}.1, namely, that the domain of definition of the various operators present in our abstract Nitsche framework can be extended to the enlarged space, $\tilde{\mathcal{V}}^{\VM} + \mathcal{V}^{\VM}_h$, where $\mathcal{V}^{\VM}_h \subset \mathcal{V}^{\VM}$. The smoothness concerns are then averted entirely because both the infinite-dimensional solution and its finite-dimensional approximation reside inside of this enlarged space.

To ground this concept in the vector Poisson problem we have considered thus far, we begin with the observation that,
\begin{equation*}
 \int_{\physBoundary_D} \left[ \left( \nabla {\bf w} \right) \cdot {\bf n} \right] \cdot {\bf v} \ d \physBoundary = \sum_{E \in \mathcal{E}_D} \int_E \left[ \left( \nabla {\bf w} \right) \cdot {\bf n} \right] \cdot {\bf v} \ d \physBoundary,
\end{equation*}
$\forall \ {\bf w} \in \tilde{\mathcal{V}}^{\VM}$ and ${\bf v} \in \mathcal{V}^{\VM}$. Therefore, we can equivalently define $\mathcal{B}^{\VM} \colon \tilde{\mathcal{V}}^{\VM} \rightarrow \left( \mathcal{Q}^{\VM} \right)^\ast$ via
\begin{align}
\langle \mathcal{B}^{\VM} {\bf w}, \mathcal{T}^{\VM} {\bf v} \rangle &\equiv \sum_{E \in \mathcal{E}_D} \int_E \left[ \left( \nabla {\bf w} \right) \cdot {\bf n} \right] \cdot {\bf v} \ d\physBoundary \nonumber
\end{align}
$\forall \ {\bf w} \in \tilde{\mathcal{V}}^{\VM}$ and ${\bf v} \in \mathcal{V}^{\VM}$. Through this alternative definition, it is clear that we can extend the domain of definition of $\mathcal{B}^{\VM}$ to functions ${\bf w} \in \mathcal{V}^{\VM}$ with square-integrable normal derivative over each edge in the Dirichlet edge mesh. In particular, we can extend it to
\begin{equation}
  \tilde{\mathcal{V}}^{\VM} + \mathcal{V}^{\VM}_h := \left\{ \tilde{\bf w} + {\bf w}_h \colon \tilde{\bf w} \in \tilde{\mathcal{V}}^{\VM} \text{ and } {\bf w}_h \in \mathcal{V}^{\VM}_h \right\},
  \label{eqn:poisson_enlarged_space}
\end{equation}
since the functions in $\mathcal{V}^{\VM}_h$ are smooth over each element in the mesh.

For the vector Poisson model problem considered in this section, the splitting of the boundary integral into a summation over the edge mesh is actually not necessary for domain extension. This is because the functions in $\mathcal{V}^{\VM}_h$ are comprised of at least $C^0$-continuous piecewise polynomials or rationals, hence, their normal derivatives are in $\left[ L^2 (\physDomain) \right]^\uDim$. Thus, the domain of definition of $\mathcal{B}^{\VM}$ trivially extends to $\tilde{\mathcal{V}}^{\VM} + \mathcal{V}^{\VM}_h$, without the need to split the boundary integral to a summation of integrals over the Dirichlet edge mesh. However, in Section~\ref{sec:biharmonic} and Section~\ref{sec:KL_plate}, we will see that this process is \emph{necessary} for the formulation to be well-defined.

For those problems, third derivatives along the boundary are required and, without such a splitting, this would render low-ordered discretizations with differential jumps at the elemental junctions inadmissible. The same approach is applied for the linear Kirchhoff-Love shell presented in \cite[\S 3.3]{Benzaken2020}. In general, for high-order PDEs the extension is necessary because the boundary integrals appearing in the Nitsche formulation may not be well defined due to the presence of high-ordered derivatives. Certain discretizations such as those containing extraordinary points or collapsed nodes resulting in degenerate finite elements may also degrade smoothness and require such a splitting. However, the elementwise-split integral is well defined for adequate polynomial or rational approximation over the mesh $\mathcal{K}$, these types of discretizations have well-defined derivatives and are $C^\infty$ away from problematic points. Henceforth, we will split all boundary integrals when presenting our Nitsche formulation, even when it is not strictly necessary for well-posedness.

\subsection{Generalized Trace and Cauchy-Schwarz Inequalities}
\label{sec:model_TI_CS}
In the third step of the construction of Nitsche's method, suitable linear maps must be specified in order to establish the generalized trace and Cauchy-Schwarz inequalities appearing in Assumption~\ref{assumption2}. To define these maps, we must introduce additional assumptions and concepts.

To begin, we consider the following lemma:

\begin{lemma}[Trace Inequality]
\label{lemma:trace_heat_prelim}
  There exists a positive dimensionless constant $C^{\VM}_{\textup{tr}} > 0$ such that
  \begin{equation}
    \sum_{E \in \mathcal{E}_D} \int_E h_E \left| \left( \nabla {\bf v}_h \right) \cdot {\bf n} \right|^2 d \physBoundary \leq C^{\VM}_{\textup{tr}} a^{\VM} ( {\bf v}_h , {\bf v}_h )
    \label{eq:Poisson_Trace}
  \end{equation}
  for all ${\bf v}_h \in \mathcal{V}^{\VM}_h$. Furthermore, an upper bound for $C^{\VM}_{\textup{tr}}$ can be obtained through an associated generalized eigenvalue problem.
  \begin{proof}

We begin by denoting the kernel of the gradient operator by
\begin{equation*}
  \text{ker}(\nabla) \equiv \left\{ {\bf v}_h \in \mathcal{V}_h^{\VM} \colon \nabla {\bf v}_h = {\bf 0} \right\}.
\end{equation*}
In this instance, the kernel of the gradient operator is simply the space of constant vectors. However in the later sections, the associated differential operators are more complex and hence, so are the kernels thereof. We furthermore denote the orthogonal complement of this space by
\begin{equation*}
  \mathring{\mathcal{V}}^{\VM}_h \equiv \left\{ {\bf v} \in \mathcal{V}_h^{\VM} \colon ({\bf v},{\bf r})_{L^2} = 0 \ \forall \ {\bf r} \in \text{ker}(\nabla) \right\}.
\end{equation*}
By construction, we are then able to express $\mathcal{V}^{\VM}_h = \mathring{\mathcal{V}}^{\VM}_h \oplus \text{ker}(\nabla)$. Consequently, for any ${\bf v}_h \in \mathcal{V}^{\VM}_h$, there exists $\mathring{\bf v}_h \in \mathring{\mathcal{V}}^{\VM}_h$ such that $\nabla {\bf v}_h = \nabla \mathring{\bf v}_h$. If there exists a positive dimensionless constant $C^{\VM}_{\textup{tr}} > 0$ such that \eqref{eq:Poisson_Trace} holds for all $\mathring{\bf v}_h \in \mathring{\mathcal{V}}^{\VM}_h$, then \eqref{eq:Poisson_Trace} therefore holds with the same constant $C^{\VM}_{\textup{tr}}$ for all ${\bf v}_h \in \mathcal{V}^{\VM}_h$.

Define the mesh-dependent bilinear form for all ${\bf w}_h, {\bf v}_h \in \mathcal{V}^{\VM}_h$ via
\begin{equation}
  b_h^{\VM}( {\bf w}_h , {\bf v}_h ) := \sum_{E \in \mathcal{E}_D} \int_E h_E \left[ \left( \nabla {\bf w}_h \right) \cdot {\bf n} \right] \cdot \left[ \left( \nabla {\bf v}_h \right) \cdot {\bf n} \right] d \physBoundary
\end{equation}
and consider the generalized eigenproblem: Find $({\bf u}_h, \lambda_h) \in \mathring{\mathcal{V}}^{\VM}_h \times \mathbb{R}$ such that
\begin{equation}
  b_h^{\VM}( {\bf u}_h , \delta {\bf u}_h ) = \lambda_h a^{\VM} ( {\bf u}_h , \delta {\bf u}_h )
  \label{eq:Gen_EP_M}
\end{equation}
for all $\delta {\bf u}_h \in \mathring{\mathcal{V}}^{\VM}_h$. By construction, the bilinear form $a^{\VM} ( \cdot , \cdot )$ is coercive on $\mathring{\mathcal{V}}^{\VM}_h$ and the left-hand side of \eqref{eq:Gen_EP_M} is bounded on $\mathring{\mathcal{V}}^{\VM}_h$. Therefore all eigenvalues of the above generalized eigenproblem are non-negative and finite. Moreover, since $\mathring{\mathcal{V}}^{\VM}_h$ is finite-dimensional, there are a finite number of eigenvalues associated with the generalized eigenproblem.
By the min-max theorem, the max eigenvalue therefore satisfies
\begin{equation}
  \lambda^{\VM}_{\textup{max}} = \sup_{\substack{{\bf v}_h \in \mathring{\mathcal{V}}^{\VM}_h \\ {\bf v}_h \neq {\bf 0}}} \frac{ b_h^{\VM}( {\bf v}_h , {\bf v}_h ) }{a^{\VM} ( {\bf v}_h, {\bf v}_h ) }.
  \label{eqn:RQ_VM}
\end{equation}
It is easily seen then that the lemma is satisfied for $C^{\VM}_{\textup{tr}} = \lambda^{\VM}_{\textup{max}}$.

  \end{proof}
\end{lemma}

\begin{remark}
  From its proof, we see that Lemma~\ref{lemma:trace_heat_prelim} is satisfied for $C^{\VM}_{\textup{tr}} = \lambda^{\VM}_{\textup{max}}$, where $\lambda^{\VM}_{\textup{max}}$ is the largest eigenvalue of the generalized eigenproblem \eqref{eq:Gen_EP_M}. However, given a basis for $\mathcal{V}^{\VM}_h$, it is in general very difficult to construct a basis for the space $\mathring{\mathcal{V}}^{\VM}_h$. If instead we consider the same eigenproblem, but posed over entirety of $\mathcal{V}^{\VM}_h$, we immediately observe that the associated bilinear form $a^{\VM}( \cdot, \cdot)$ is not coercive over this space. Therefore $\lambda^{\VM}_{\textup{max}} = \infty$ by \eqref{eqn:RQ_VM}, and the corresponding eigenvector is in $\text{ker}(\nabla)$. Fortunately, recalling the orthogonal decomposition $\mathcal{V}^{\VM}_h = \mathring{\mathcal{V}}^{\VM}_h \oplus \text{ker}(\nabla)$ presented in the proof of Lemma~\ref{lemma:trace_heat_prelim}, the largest \emph{finite} eigenvalue of the eigenprobem posed over $\mathcal{V}^{\VM}_h$ corresponds to the largest eigenvalue of the same eigenproblem but posed over $\mathring{\mathcal{V}}^{\VM}_h$. This observation can be used to provide an explicit estimate of $C^{\VM}_{\textup{tr}}$ without needing to construct a basis for $\mathring{\mathcal{V}}^{\VM}_h$.

 In particular, given a basis $\{ {\bf N}_i \}_{i=1}^n$ for the space $\mathcal{V}^{\VM}_h$, it follows that $\lambda^{\VM}_{\textup{max}}$ may be computed as the largest finite eigenvalue of the generalized matrix eigenproblem $\left({\bf B} - \lambda_h {\bf A} \right) {\bf x} = {\bf 0}$, where
 \begin{equation*}
  \left[ {\bf A} \right]_{ij} := a^{\VM} ( {\bf N}_i, {\bf N}_j ), \hspace{15pt} \left[ {\bf B} \right]_{ij} := b_h^{\VM} ( {\bf N}_i , {\bf N}_j ),
\end{equation*}
 and ${\bf x} = \{ x_i \}_i$ is the set of coefficients such that ${\bf u}_h = \sum_i x_i {\bf N}_i$ is an eigenvector with associated eigenvalue $\lambda_h$.Note that since this matrix eigenvalue problem has infinite eigenvalues, special care may be needed for their explicit numerical computation.
  \label{remark:model_EP_for_TI}
\end{remark}

With the above assumptions and concepts in hand, we are now ready to define suitable linear maps such that the generalized trace and Cauchy-Schwarz inequalities given in Assumption~\ref{assumption2} are satisfied.

Beginning with the concept of the generalized trace inequality, we first define $\eta^{\VM}: \textup{dom}(\eta^{\VM}) \subseteq \left( \mathcal{Q}^{\VM} \right)^{\ast} \rightarrow \mathcal{Q}^{\VM}$ to be a densely defined, positive, self-adjoint linear map on the space
\begin{equation*}
\left\{ \mathcal{B}^{\VM} {\bf v}: {\bf v} \in \tilde{\mathcal{V}}^{\VM} + \mathcal{V}^{\VM}_h \right\}
\end{equation*}
that satisfies
\begin{equation}
\left\langle \mathcal{B}^{\VM} {\bf w}, \eta^{\VM} \mathcal{B}^{\VM} {\bf v} \right\rangle \equiv \sum_{E \in \mathcal{E}_D} \int_E \frac{h_E }{C^{\VM}_{\textup{tr}}} \left[ \left( \nabla {\bf w}_h \right) \cdot {\bf n} \right] \cdot \left[ \left( \nabla {\bf v}_h \right) \cdot {\bf n} \right] d \physBoundary
\label{eqn:etaDefM}
\end{equation}
for all ${\bf w}, {\bf v} \in \tilde{\mathcal{V}}^{\VM} + \mathcal{V}^{\VM}_h$. This implicit definition of $\eta^{\VM}$ is inspired by the trace inequality problem as presented in Lemma~\ref{lemma:trace_heat_prelim}, which leads to the following lemma.

\begin{lemma}[Generalized Trace Inequality for the Poisson Problem]
\label{lemma:trace_heat}
  It holds that
  \begin{equation*}
  \langle \mathcal{B}^{\VM} {\bf v}_h, \eta^{\VM} \mathcal{B}^{\VM} {\bf v}_h \rangle \leq a^{\VM}({\bf v}_h,{\bf v}_h)
\end{equation*}
  for all ${\bf v}_h \in \mathcal{V}^{\VM}_h$.
  \begin{proof}
  The proof immediately follows from Lemma~\ref{lemma:trace_heat_prelim} and the definition of $\eta^{\VM}$.
  \end{proof}
\end{lemma}

Next, we define the linear map $\epsilon^{\VM}: \textup{dom}(\epsilon^{\VM}) \subseteq \left( \mathcal{Q}^{\VM} \right)^{\ast} \rightarrow \mathcal{Q}^{\VM}$ through the action of its inverse as
\begin{equation}
\left\langle \left( \epsilon^{\VM} \right)^{-1} {\bf w}, {\bf v} \right\rangle := \sum_{E \in \mathcal{E}_D} \int_E \frac{C^{\VM}_\text{pen}}{h_E} {\bf w} \cdot {\bf v} \ d \physBoundary
\label{eqn:epsDefVM}
\end{equation}
for all ${\bf w}, {\bf v} \in \mathcal{Q}^{\VM}$, where $C^{\VM}_\text{pen} > C^{\VM}_\text{tr}$ is a chosen dimensionless constant. Then we have the following result:

\begin{lemma}[Generalized Cauchy-Schwarz Inequality for the Poisson Problem]
\label{lemma:cs_heat}
  Let $C^{\VM}_{\textup{pen}} = \CSconstant^2 C^{\VM}_{\textup{tr}}$, where $\CSconstant \in (1,\infty)$.  Then
  \begin{equation*}
  \left| \left\langle \mathcal{B}^{\VM} {\bf v}, \mathcal{T}^{\VM} {\bf w} \right\rangle \right| \leq \frac{1}{\CSconstant} \left\langle \mathcal{B}^{\VM} {\bf v}, \eta^{\VM} \mathcal{B}^{\VM} {\bf v} \right\rangle^{1/2} \left\langle \left(\epsilon^{\VM}\right)^{-1} \mathcal{T}^{\VM} {\bf w}, \mathcal{T}^{\VM} {\bf w} \right\rangle^{1/2} \label{eq:Poisson_CS}
\end{equation*}
  for all ${\bf v}, {\bf w} \in \tilde{\mathcal{V}}^{\VM} + \mathcal{V}^{\VM}_h$.
  \begin{proof}
  We write
  \begin{align*}
  \left\langle \mathcal{B}^{\VM} {\bf v}, \mathcal{T}^{\VM} {\bf w} \right\rangle &= \sum_{E \in \mathcal{E}_D} \int_E \left[ \left( \nabla {\bf v} \right) \cdot {\bf n} \right] \cdot {\bf w} \ d\physBoundary \nonumber \\
  &= \sum_{E \in \mathcal{E}_D} \int_E \frac{1}{\CSconstant^2} \frac{h_E}{\cTrace{}^{\VM}} \frac{\cPen{}^{\VM}}{h_E} \left[ \left( \nabla {\bf v} \right) \cdot {\bf n} \right] \cdot {\bf w} \ d\physBoundary \nonumber \\
  &\leq \frac{1}{\CSconstant} \sum_{E \in \mathcal{E}_D} \left( \int_E \frac{h_E}{C^{\VM}_{\textup{tr}} } \left| \left( \nabla {\bf v} \right) \cdot {\bf n} \right|^2 d\physBoundary \right)^{1/2} \left( \int_E \frac{C^{\VM}_{\textup{pen}}}{h_E} \left| {\bf w} \right|^2 d\physBoundary \right)^{1/2} \nonumber \\
    &\leq \frac{1}{\CSconstant} \left( \sum_{E \in \mathcal{E}_D} \int_E \frac{h_E}{C^{\VM}_{\textup{tr}}} \left| \left( \nabla {\bf v} \right) \cdot {\bf n} \right|^2 d\physBoundary \right)^{1/2} \left( \sum_{E \in \mathcal{E}_D} \int_E \frac{C^{\VM}_{\textup{pen}}}{h_E} \left| {\bf w} \right|^2 d\physBoundary \right)^{1/2} \nonumber \\
    &= \frac{1}{\CSconstant} \left\langle \mathcal{B}^{\VM} {\bf v}, \eta^{\VM} \mathcal{B}^{\VM} {\bf v} \right\rangle^{1/2} \left\langle \left(\epsilon^{\VM}\right)^{-1} \mathcal{T}^{\VM} {\bf w}, \mathcal{T}^{\VM} {\bf w} \right\rangle^{1/2},
  \end{align*}
    where the standard continuous Cauchy-Schwarz inequality was employed in the first inequality above ($(f, g)_{L^2(D)} \leq \| f \|_{L^2(D)} \| g \|_{L^2(D)}$ for $f, g \in L^2(D)$) and the standard discrete Cauchy-Schwarz inequality was employed in the second inequality above ($|(x,y)|_{\ell^2} \leq \| x \|_{\ell^2} \| y \|_{\ell^2}$ for $x, y \in \mathbb{R}^N$ for $N \in \N$).
  \end{proof}
\end{lemma}

With these selections for $\eta^{\VM}$ and $\epsilon^{\VM}$, the generalized trace and Cauchy-Schwarz inequalities appearing in Assumption~\ref{assumption2} are satisfied.

\begin{remark}
According to our analysis, a practitioner may select $C^{\VM}_{\textup{pen}} = \CSconstant^2 C^{\VM}_{\textup{tr}}$ for any $\CSconstant \in (1,\infty)$.  Generally speaking, the Dirichlet boundary condition is enforced more strongly for larger $\CSconstant$ as opposed to smaller $\CSconstant$.  However, the condition number of the linear system associated with Nitsche's method scales linearly with $\CSconstant$ \cite{juntunen2009nitsche} and, in certain circumstances, the discrete solution becomes over-constrained and boundary locking occurs as $\CSconstant \rightarrow \infty$, resulting in a loss of solution accuracy \cite{lew2008discontinuous}.  On the other hand, as $\CSconstant \rightarrow 1$, the linear system associated with Nitsche's method may lose definiteness, and Theorem~\ref{theorem:convergence_heat} suggests that the energy norm error may blow up in the limit $\CSconstant \rightarrow 1$.  It is therefore advisable to choose a moderate value for $\CSconstant$.  Based on our collective experience, we recommend setting $\CSconstant = 2$.
\end{remark}

\begin{remark}
  The definitions of the linear maps $\eta^{\VM}$ and $\epsilon^{\VM}$ are simple for the Poisson problem. In fact, their actions can be entirely represented in this case through the elementwise definitions
  \begin{equation}
    \eta^{\VM}_E \equiv \frac{h_E}{C^{\VM}_\text{tr}} \hspace{15pt} \text{and} \hspace{15pt} \epsilon^{\VM}_E \equiv \frac{h_E}{C^{\VM}_\text{pen}}.
  \end{equation}
  Then the linear maps $\eta^{\VM}$ and $\epsilon^{\VM}$ are defined simply as the summation over these maps as seen in \eqref{eqn:etaDefM} and \eqref{eqn:epsDefVM}, respectively. However, in general, these linear maps cannot be represented so concisely as we will see in Section~\ref{sec:TI_CS_B}.
  \label{rem:eps_eta_VM}
\end{remark}

\subsection{Nitsche's Method}

In the third and final step of the construction of Nitsche's method for a given variational problem, we simply pose Nitsche's method according to $(N_h)$.

\vspace{1em}
\begin{mybox}[\emph{Nitsche's Method for the Poisson Problem}]
  \vspace{-3pt}
$$
(N_h^{\VM}) \left\{ \hspace{5pt}
\parbox{4.35in}{
Given $\linFunctional^{\VM} \in \left( \mathcal{V}^{\VM} \right)^*$ and ${\bf g} \in \mathcal{Q}^{\VM}$, find ${\bf u}_h \in \mathcal{V}^{\VM}_h$ such that
\begin{equation*}
\begin{aligned}
a_h^{\VM}({\bf u}_h,\delta {\bf u}_h) &= \underbrace{ \int_{\physDomain} {\bf f} \cdot \delta {\bf u}_h \ d \physDomain + \int_{\physBoundary_N} {\bf h} \cdot \delta {\bf u}_h \ d \physBoundary }_{ \displaystyle \langle \linFunctional^{\VM}, \delta {\bf u}_h \rangle}\\
&\phantom{=} + \sum_{E \in \mathcal{E}_D} \left( {\color{ForestGreen} \underbrace{ - \int_E \left[ \left( \nabla \delta {\bf u}_h \right) \cdot {\bf n} \right] \cdot {\bf g}  \ d \physBoundary }_\text{Symmetry Term} } {\color{Orchid} \underbrace{ + \int_E \frac{\cPen{}^{\VM}}{h_E} {\bf g} \cdot \delta {\bf u}_h \ d\Gamma}_\text{Penalty Term} } \right)
\end{aligned}
\label{eqn:Linear_Elasticity_Weak_Nitsche}
\end{equation*}
for every $\delta {\bf u}_h \in \mathcal{V}^{\VM}_h$ where $a^{\VM}_h: \left( \tilde{\mathcal{V}}^{\VM} + \mathcal{V}^{\VM}_h \right) \times \left( \tilde{\mathcal{V}}^{\VM} + \mathcal{V}^{\VM}_h \right) \rightarrow \mathbb{R}$ is the bilinear form defined by
\begin{equation*}
\begin{aligned}
a^{\VM}_h({\bf w}_h, {\bf v}_h) &\equiv \underbrace{ \int_{\physDomain} \left( \nabla {\bf w}_h \right) : \left( \nabla {\bf v}_h \right) \ d \physDomain }_{ \displaystyle a^{\VM}({\bf w}_h, {\bf v}_h) } + \sum_{E \in \mathcal{E}_D} \left( {\color{Cerulean} \underbrace{ - \int_E \left[ \left( \nabla {\bf w}_h \right) \cdot {\bf n} \right] \cdot {\bf v}_h \ d \physBoundary }_\text{Consistency Term} }  \right.\\
&\phantom{\equiv} \left. {\color{ForestGreen} \underbrace{ - \int_E \left[ \left( \nabla {\bf v}_h \right) \cdot {\bf n} \right] \cdot {\bf w}_h \ d \physBoundary }_\text{Symmetry Term} } {\color{Orchid} \underbrace{ + \int_E \frac{\cPen{}^{\VM}}{h_E} {\bf w}_h \cdot {\bf v}_h \ d\Gamma}_\text{Penalty Term} } \right).
\end{aligned}
\end{equation*}
for ${\bf w}_h, {\bf v}_h \in \tilde{\mathcal{V}}^{\VM} + \mathcal{V}^{\VM}_h$.
}
\right.
\vspace{3pt}
$$
\end{mybox}
\vspace{1em}

\noindent Since Assumptions~\ref{assumption1} and~\ref{assumption2} from the abstract framework in Section~\ref{sec:Nitsche} are satisfied as a consequence of Lemmas~\ref{lemma:greens_heat},~\ref{lemma:trace_heat}, and~\ref{lemma:cs_heat}, we have the following theorem stating well-posedness and an error estimate for Nitsche's method for the Poisson problem:

\begin{theorem}[Well-Posedness and Error Estimate for the Poisson Problem]
\label{theorem:convergence_heat}
  Let $C^{\VM}_{\textup{pen}} = \CSconstant^2 C^{\VM}_{\textup{tr}}$, where $\CSconstant \in (1,\infty)$.  Then there exists a unique discrete solution ${\bf u}_h \in \mathcal{V}^{\VM}_h$ to Problem $(N_h^{\VM})$.  Moreover, if the continuous solution ${\bf u} \in \mathcal{V}^{\VM}$ to Problem $(V^{\VM})$ satisfies ${\bf u} \in \tilde{\mathcal{V}}^{\VM}$, then the discrete solution ${\bf u}_h$ satisfies the error estimate
  \begin{equation*}
\vvvertiii{{\bf u} - {\bf u}_h}_{\VM} \leq \left( 1+ \frac{2}{1-\frac{1}{\CSconstant}} \right) \min_{{\bf v}_h \in \mathcal{V}^{\VM}_h} \vvvertiii{{\bf u} - {\bf v}_h}_{\VM}.
 \label{eq:Poisson_Error}
  \end{equation*}
  where
  \begin{equation}
    \vvvertiii{{\bf w}}^2_{\VM} \equiv \int_{\physDomain} \left| \nabla {\bf w} \right|^2 d \physDomain + \sum_{E \in \mathcal{E}_D} \left( \int_E \frac{h_E}{\cTrace{}^{\VM}} \left| \left( \nabla {\bf w} \right) \cdot {\bf n} \right|^2 d \physBoundary + 2 \int_E \frac{\cPen{}^{\VM}}{h_E} \left| {\bf w} \right|^2 d\Gamma \right).
  \end{equation}

  \begin{proof}
  Note that Nitsche's method for the Poisson problem precisely fits into the abstract variational framework presented in Section~\ref{sec:Nitsche} with $\mathcal{V} = \mathcal{V}^{\VM}$, $\mathcal{Q} = \mathcal{Q}^{\VM}$, $a(\cdot,\cdot) = a^{\VM}(\cdot,\cdot)$, $\linFunctional = \linFunctional^{\VM}$, $\mathcal{T} = \mathcal{T}^{\VM}$, $\tilde{\mathcal{V}} = \tilde{\mathcal{V}}^{\VM}$, $\mathcal{L} = \mathcal{L}^{\VM}$, $\mathcal{B} = \mathcal{B}^{\VM}$, $\mathcal{V}_h = \mathcal{V}^{\VM}_h$, $\epsilon = \epsilon^{\VM}$, and $\eta = \eta^{\VM}$.  Moreover, Assumption~\ref{assumption1} of the abstract variational framework is satisfied due to Lemma~\ref{lemma:greens_heat}, and Assumption~\ref{assumption2} is satisfied due to Lemmas~\ref{lemma:trace_heat} and~\ref{lemma:cs_heat}.  Well-posedness and the error estimate therefore follow directly from Theorem~\ref{theorem:error_estimate}.
  \end{proof}
\end{theorem}

The above result indicates that Nitsche's method is quasi-optimal in the energy norm (in the sense that the error in the discrete solution is proportional to the best approximation error) when the continuous solution ${\bf u} \in \mathcal{V}^{\VM}$ to Problem $(V^{\VM})$ satisfies ${\bf u} \in \tilde{\mathcal{V}}^{\VM}$.  However, the above result does not reveal the rates of convergence of the energy norm error for Nitsche's method, nor does it reveal the rates of convergence for other norms one may care about (for instance, the $L^2$-norm).  Such a convergence analysis can be quite technical for complex problems and requires the use of both interpolation and trace estimates, which is beyond the scope of this paper.

\begin{remark}
  The presented formulation of Nitsche's method for the Poisson problem is standard, even though we derived it from our abstract framework.  Well-posedness and error estimates for the formulation have also been proved previously using similar techniques to those presented here \cite{stenberg1995some}.  The advantage of the abstract framework is that it streamlines the construction and mathematical analysis of Nitsche's method for more complicated problems.
\end{remark}


\section{Nitsche's Method for the Biharmonic Equation}
\label{sec:biharmonic}

Having completed the rigorous derivation of Nitsche's method for the Poisson problem, we consider the more complicated vector-valued biharmonic problem in this section. As in the previous section, we utilize our abstract framework and show that it is relatively straight-forward to apply it to high-order PDEs. For the biharmonic problem  considered here, the $4^{th}$-order nature presents two sets of boundary conditions that are easily handled as long as each step identified at the end of Section~\ref{sec:Nitsche} is performed carefully. The steps necessary to obtain the Green's identity are more complex, but still follow a similar pattern to that of the Green's identity for the Poisson problem. We also observe that the linear maps $\epsilon$ and $\eta$ in this case can only be concisely written in terms of their action. This is in contrast to the result in the previous section where these linear maps were expressed explicitly with relative ease.


\subsection{The Variational Formulation}

Starting with the first step in our recipe once again, we begin by stating the minimization problem and its corresponding variational formulation. Proceeding forward, we use superscript $\B$ to refer to the \emph{biharmonic} problem considered here.

Let $\physDomain \subset \mathbb{R}^\sDim$ be a bounded domain with Lipschitz-continuous boundary $\physBoundary = \partial \physDomain$, where $\sDim \in \N$ is the spatial dimension.  We further assume that $\physDomain$ is $H^2$-regular (i.e., we assume that $H^2(\Omega) \equiv \left\{ v \in H^1(\Omega) : \Delta v \in L^2(\Omega) \right\}$). This precludes, for example, domains with re-entrant corners. Given this, let ${\bf u}$, ${\bf v}$, and ${\bf w}$ denote vector-valued functions of dimension $\uDim \in \N$. We are interested in minimizing a functional with an internal energy given by
\begin{eqnarray*}
E^{\B}_{\textup{int}}({\bf u}) \equiv \frac{1}{2} \int_{\physDomain} \left( \Delta {\bf u} \right) \cdot \left( \Delta {\bf u} \right) \ d\physDomain.
\end{eqnarray*}
First, we must select a space of admissible solutions that is smooth enough such that this functional is well-defined. We therefore restrict our attention to functions with at least \emph{two} integrable derivatives. Thus, let $\mathcal{V}^{\B} \equiv \left[ H^2(\physDomain) \right]^\uDim$ be the space of admissible solutions free of boundary conditions.

Next, we partition the boundary into Dirichlet and Neumann parts. In particular, let $\physBoundary_{D_1}$ and $\physBoundary_{N_1}$ denote what we refer to as the Dirichlet-1 and Neumann-1 boundaries associated with the prescribed function values and the normal derivative of the Laplacian, respectively. Similarly, let $\physBoundary_{D_2}$ and $\physBoundary_{N_2}$ denote what we refer to as the Dirichlet-2 and Neumann-2 boundaries associated with prescribed normal derivatives and Laplacians, respectively. For $\alpha = 1,2$, we require that $\physBoundary = \overline{\physBoundary_{D_\alpha} \cup \physBoundary_{N_\alpha}}$, $\physBoundary_{D_\alpha} \cap \physBoundary_{N_\alpha} = \emptyset$, and $\physBoundary_{D_\alpha} \neq \emptyset$ for a well-posed PDE. Note that there are no constraints between the 1- and 2-boundaries because there is no energetic exchange between the two sets, that is, the function value and derivative at the boundary can be prescribed entirely independent of each other.

With the above boundary partitions in place and the selection for $\mathcal{V}^\B$, the associated trace space is then $\mathcal{Q}^{\B} \equiv \left[ H^{3/2}(\physBoundary_{D_1})\right]^\uDim \times \left[ H^{1/2}(\physBoundary_{D_2})\right]^\uDim$. Furthermore, we can now define the trace operator $\mathcal{T}^{\B} \colon \mathcal{V}^{\B} \rightarrow \mathcal{Q}^{\B}$ via its action on an arbitrary ${\bf w} \in \mathcal{V}^{\B}$ via $\mathcal{T}^{\B} {\bf w} \equiv \left( {\bf w} \big|_{\physBoundary_{D_1}}, \nabla {\bf w} \cdot {\bf n} \big|_{\physBoundary_{D_2}} \right)$ or, in other words, ${\bf w}$, and its normal derivative, $\nabla {\bf w} \cdot {\bf n}$, restricted to their Dirichlet boundaries. In light of this, we are able to prescribe Dirichlet boundary conditions provided they are smooth enough. Given a prescribed displacement and normal derivative on the Dirichlet boundaries, $(\textup{\bf g}, \textup{\bf h}) \in \mathcal{Q}^{\B}$, we denote the trial space of displacement fields satisfying the prescribed Dirichlet boundary conditions by
\begin{equation*}
\mathcal{V}^{\B}_{\textup{\bf g},\textup{\bf h}} \equiv \left\{ {\bf v} \in \mathcal{V}^{\B} \colon \mathcal{T}^{\B} {\bf v} = (\textup{\bf g},\textup{\bf h}) \right\},
\end{equation*}
Similarly, $\mathcal{V}^{\B}_{{\bf 0},{\bf 0}}$ denotes the test space of displacement fields satisfying homogeneous Dirichlet boundary conditions, corresponding to $\textup{\bf g} \equiv {\bf 0}$ and $\textup{\bf h} \equiv {\bf 0}$.

With only homogeneous forcing functions, the minimizer of $E^{\B}_{\textup{int}}(\cdot)$ is dictated purely by the choice of the Dirichlet boundary conditions, ${\bf g}$ and ${\bf h}$. To accommodate the additional scenarios that comprise the general biharmonic problem, we additionally consider the functional representing the external work, which takes the form
\begin{eqnarray*}
E^{\B}_{\textup{ext}}({\bf u}) \equiv - \int_{\physDomain} \textup{\bf f} \cdot {\bf u} \ d \physDomain + \int_{\physBoundary_{N_1}} \textup{\bf p} \cdot {\bf u} \ d \physBoundary - \int_{\physBoundary_{N_2}} \textup{\bf q} \cdot \left[ \left( \nabla {\bf u} \right) \cdot {\bf n} \right] \ d \physBoundary.
\end{eqnarray*}
For this functional to be well-defined, we require that $\textup{\bf f} \in \left[ L^2(\physDomain) \right]^\uDim$, $\textup{\bf p} \in \left[ L^2(\physBoundary_{N_1}) \right]^\uDim$, and $\textup{\bf q} \in \left[ L^2(\physBoundary_{N_2}) \right]^\uDim$. The sign convention for the terms in $E^\B_{\textup{ext}}({\bf u})$ will become apparent later when we present the generalized Green's identity associated with this problem.

We define an associated bilinear form $a^{\B}(\cdot,\cdot): \mathcal{V}^{\B} \times \mathcal{V}^{\B} \rightarrow \mathbb{R}$ via
\begin{equation*}
a^{\B}({\bf w},{\bf v}) \equiv \int_{\physDomain} \left( \Delta {\bf w} \right) \cdot \left( \Delta {\bf v} \right) \ d\physDomain
\label{eqn:aB}
\end{equation*}
for all ${\bf w},{\bf v} \in \mathcal{V}^{\B}$. The linear functional $\linFunctional^{\B} \in \left( \mathcal{V}^{\B} \right)^{\ast}$ is defined as
\begin{equation*}
\left\langle \linFunctional^{\B}, {\bf v} \right\rangle \equiv -  E^{\B}_{\textup{ext}}({\bf v})
\label{eqn:fB}
\end{equation*}
for all ${\bf v} \in \mathcal{V}^{\B}$.

Given all of these definitions, the minimization problem is simply
$$
(M^{\B}) \left\{ \hspace{5pt}
\parbox{4.35in}{
\noindent Find ${\bf u} \in \mathcal{V}^{\B}_{\textup{\bf g},\textup{\bf h}}$ that minimizes the total energy
\begin{eqnarray*}
E^{\B}_{\textup{total}}({\bf u}) = E^{\B}_{\textup{int}}({\bf u}) + E^{\B}_{\textup{ext}}({\bf u}).
\end{eqnarray*}
}
\right.
$$

The solution to Problem $(M^{\B})$ is also the solution to the following variational problem:
$$
(V^{\B}) \left\{ \hspace{5pt}
\parbox{4.35in}{
\noindent Find ${\bf u} \in \mathcal{V}^{\B}_{\textup{\bf g},\textup{\bf h}}$ such that
\begin{eqnarray*}
a^{\B}({\bf u}, \delta {\bf u}) = \left\langle \linFunctional^{\B}, \delta {\bf u} \right\rangle \label{eq:biharmonic_weak}
\end{eqnarray*}
for every $\delta {\bf u} \in \mathcal{V}^{\B}_{{\bf 0},{\bf 0}}$.
}
\right.
$$

\subsection{A Generalized Green's Identity}

The second step in our recipe for constructing a Nitsche-based formulation requires a generalized Green's identity that satisfies Assumption~\ref{assumption1} to be established. Similar to Section~\ref{sec:greenpoisson}, we start with a product rule, but this time we invoke it twice to arrive at the relevant identity. First, observe that
\begin{equation*}
    \nabla \cdot \left[ \left( \nabla \Delta {\bf w} \right) \cdot {\bf v} \right] = \left[ \nabla \cdot \left( \nabla \Delta {\bf w} \right) \right] \cdot {\bf v} + \left( \nabla \Delta {\bf w} \right) : \left( \nabla {\bf v} \right)
\end{equation*}
and
\begin{equation*}
    \nabla \cdot \left[ \left( \Delta {\bf w} \right) \cdot \left( \nabla {\bf v} \right) \right] = \left( \nabla \Delta {\bf w} \right) : \left( \nabla {\bf v} \right) + \left( \Delta {\bf w} \right) \cdot \left( \Delta {\bf v} \right).
\end{equation*}
Solving the second expression for $\left( \nabla \Delta {\bf w} \right) : \left( \nabla {\bf v} \right)$ and inserting it into the first expression yields the following identity:
\begin{equation}
    \nabla \cdot \left[ \left( \nabla \Delta {\bf w} \right) \cdot {\bf v} \right] = \left( \Delta^2 {\bf w} \right) \cdot {\bf v} + \nabla \cdot \left[ \left( \Delta {\bf w} \right) \cdot \left( \nabla {\bf v} \right) \right] - \left( \Delta {\bf w} \right) \cdot \left( \Delta {\bf v} \right),
\label{eqn:div_prod_rule_biharmonic}
\end{equation}
where we have used the relationship $\nabla \cdot \nabla = \Delta$ to arrive at $\Delta^2 {\bf w} \equiv \Delta \Delta {\bf w}$, the biharmonic operator. Rearranging this expression, integrating both sides over the domain $\physDomain$, and applying the divergence theorem twice yields
\begin{equation}
  \begin{aligned}
    a^{\B}({\bf w},{\bf v}) &= \int_\physDomain \left( \Delta {\bf w} \right) \cdot \left( \Delta {\bf v} \right) \ d \physDomain \\
    &= \int_\physDomain \left( \Delta^2 {\bf w} \right) \cdot {\bf v} \ d \physDomain - \int_\physDomain \nabla \cdot \left[ \left( \nabla \Delta {\bf w} \right) \cdot {\bf v} \right] \ d \physDomain + \int_\physDomain \nabla \cdot \left[ \left( \Delta {\bf w} \right) \cdot \left( \nabla {\bf v} \right) \right] \ d \physDomain \\
    &= \int_\physDomain \left( \Delta^2 {\bf w} \right) \cdot {\bf v} \ d \physDomain - \int_{\physBoundary} \left[ \left( \nabla \Delta {\bf w} \right) \cdot {\bf n} \right] \cdot {\bf v} \ d \physBoundary + \int_{\physBoundary} \left( \Delta {\bf w} \right) \cdot \left[ \left( \nabla {\bf v} \right) \cdot {\bf n} \right] \ d \physBoundary,
  \end{aligned}
  \label{eqn:Greens_ID_biharmonic_derive}
\end{equation}
where ${\bf n}$ is a unit normal vector along $\physBoundary$ pointing outward from $\physDomain$. Note that although we have successfully applied the product rule to obtain a generalized Green's identity for the problems considered in this paper, it is not guaranteed to work for general PDEs. Instead, arriving at the generalized Green's identity may require a clever application of identities from vector calculus as well as identities emerging from other fields. For example, techniques from differential geometry were used in obtaining a generalized Green's identity for the linearized Kirchhoff-Love shell in \cite{Benzaken2020}.

Similar to Section~\ref{sec:greenpoisson}, we note that the last equality in \eqref{eqn:Greens_ID_biharmonic_derive} is no longer valid for general ${\bf w} \in \mathcal{V}^{\B}$. This is because for ${\bf w} \in \mathcal{V}^{\B}$, it follows that ${\bf w} \in \left[ L^2(\physBoundary_1) \right]^\uDim$ and $\left( \nabla {\bf w} \right) \cdot {\bf n} \in \left[ L^2(\physBoundary_2) \right]^\uDim$, but for \eqref{eqn:Greens_ID_biharmonic_derive} to hold, we also require that (i) $\Delta^2 {\bf w} \in \left[ L^2(\physDomain) \right]^\uDim$, (ii) $\left( \nabla \Delta {\bf w} \right) \cdot {\bf n} \in \left[ L^2(\physBoundary_1) \right]^\uDim$, and (iii) $\Delta {\bf w} \in \left[ L^2(\physBoundary_2) \right]^\uDim$. To satisfy these additional smoothness requirements, we let
\begin{equation*}
\begin{aligned}
&\tilde{\mathcal{V}}^{\B} \equiv \\
&\left\{ {\bf v} \in \left[ H^2(\physDomain) \right]^\uDim : \Delta^2 {\bf w} \in \left[ L^2(\physDomain) \right]^\uDim, \left( \nabla \Delta {\bf w} \right) \cdot {\bf n} \in \left[ L^2(\physBoundary_1) \right]^\uDim, \text{ and } \Delta {\bf w} \in \left[ L^2(\physBoundary_2) \right]^\uDim \right\}
\end{aligned}
\end{equation*}
and note that \eqref{eqn:Greens_ID_biharmonic_derive} is now valid for general ${\bf w} \in \tilde{\mathcal{V}}^{\B}$. Additionally, $\tilde{\mathcal{V}}^{\B} \subset \mathcal{V}^{\B}$, a prerequisite for Assumption~\ref{assumption1}.

Once again, we typically need to demonstrate that $\mathcal{L}^{\B} {\bf u} = \linFunctional^{\B}$ when ${\bf u} \in \tilde{\mathcal{V}}^{\B}$ to satisfy the requisites of Assumption~\ref{assumption1}. However, we omit the details and simply assert that it is true, as discussed in Remark~\ref{rem:skip_Lu_is_f}.

We are now ready to state the following lemma regarding the generalized Green's identity for the vector-valued biharmonic problem.

\begin{lemma}[Generalized Green's Identity for the Biharmonic Problem]
\label{lemma:greens_biharmonic}
  For $w \in \tilde{\mathcal{V}}^{\B}$ and $v \in \mathcal{V}^{\B}$, the following generalized Green's identity holds:
  \begin{equation}
    \begin{aligned}
    a^{\B}({\bf w},{\bf v}) &= \underbrace{ \int_{\physDomain} {\bf v} \cdot \left( \Delta^2 {\bf w} \right) \ d \physDomain - \int_{\physBoundary_{N_1}} \left[ \left( \nabla \Delta {\bf w} \right) \cdot {\bf n} \right] \cdot {\bf v} \ d \physBoundary + \int_{\physBoundary_{N_2}} \left( \Delta {\bf w} \right) \cdot \left[ \left( \nabla {\bf v} \right) \cdot {\bf n} \right] \ d \physBoundary }_{ \displaystyle \langle \mathcal{L}^{\B} {\bf w}, {\bf v} \rangle }\\
    &\phantom{=} \underbrace{ - \int_{\physBoundary_{D_1}} \left[ \left( \nabla \Delta {\bf w} \right) \cdot {\bf n} \right] \cdot {\bf v} \ d \physBoundary + \int_{\physBoundary_{D_2}} \left( \Delta {\bf w} \right) \cdot \left[ \left( \nabla {\bf v} \right) \cdot {\bf n} \right] \ d \physBoundary }_{ \displaystyle \langle \mathcal{B}^{\B} {\bf w}, \mathcal{T}^{\B} {\bf v} \rangle }.
  \end{aligned}
  \label{eqn:Greens_ID_B}
\end{equation}
Moreover, the solution ${\bf u}$ of Problem $(V^{\B})$ satisfies $\mathcal{L}^{\B} {\bf u} = \linFunctional^{\B}$ provided the problem parameters ${\bf f}$, ${\bf g}$, ${\bf h}$, ${\bf p}$, and ${\bf q}$ are smooth enough such that ${\bf u} \in \tilde{\mathcal{V}}^{\B}$.

  \begin{proof}
    The generalized Green's identity follows from the steps outlined in \eqref{eqn:Greens_ID_biharmonic_derive}. The result that $\mathcal{L}^{\B} {\bf u} = \linFunctional^{\B}$ follows by the same techniques used in \cite{Benzaken2020} as discussed in Remark~\ref{rem:skip_Lu_is_f}.
    
  \end{proof}
  \end{lemma}

  \begin{remark}
    The strong form of Problem $(V^{\B})$ is given by:
    $$
    (S^{\B}) \left\{ \hspace{5pt}
    \parbox{4.35in}{
    \noindent \textup{Find ${\bf u}: \overline{\physDomain} \rightarrow \R$ such that:}
    \begin{equation*}
    \begin{aligned}\begin{array}{rll}
    \Delta^2 {\bf u} &= \textup{\bf f} \hspace{10pt} &\textup{in} \ \physDomain\\
    {\bf u} &= \textup{\bf g} \hspace{10pt} &\textup{on} \ \physBoundary_{D_1}\\
    \left( \nabla {\bf u} \right) \cdot {\bf n} &= \textup{\bf h} \hspace{10pt} &\textup{on} \ \physBoundary_{D_2}\\
    \left( \nabla \Delta {\bf u} \right) \cdot {\bf n} &= \textup{\bf p} \hspace{10pt} &\textup{on} \ \physBoundary_{N_1}\\
    \Delta {\bf u} &= \textup{\bf q} \hspace{10pt} &\textup{on} \ \physBoundary_{N_2}\\
    \end{array}
    \end{aligned}
    \label{eqn:Biharmonic_Strong}
    \end{equation*}
    }
    \right.
    $$
    This result follows immediately from the relationship $\mathcal{L}^{\B} {\bf u} = \linFunctional^{\B}$ in Lemma~\ref{lemma:greens_biharmonic} and the essential boundary conditions encoded by the boundary operator $\mathcal{B}^{\B}$ and the space $\mathcal{V}^{\B}_{\bf g}$.
  \end{remark}

\subsection{Generalized Trace and Cauchy-Schwarz Inequalities}
\label{sec:TI_CS_B}

With a Green's identity in place, we are ready to provide generalized trace and Cauchy-Schwarz inequalities satisfying Assumption~\ref{assumption2}. Following an analogous procedure to that of our scalar and vector Poisson model problems, we establish a mesh $\mathcal{K}$ associated with $\physDomain$ that is comprised of elements such that $\physDomain = \text{int}(\overline{\cup_{K \in \mathcal{K}} K})$. Next, we assume that the approximation space $\mathcal{V}^{\B}_h$ consists of (at least) $C^1$-continuous piecewise polynomial or rational approximations over the mesh $\mathcal{K}$. We collect the boundary edges into an edge mesh $\mathcal{E}$. For the biharmonic problem, we must construct two additional edge meshes, $\mathcal{E}_{D_1}$ and $\mathcal{E}_{D_2}$. The construction of these sets follows analogously to the previous sections with the exception that we associate the members of $\mathcal{E}_{D_1}$ with elements whose edges belong to $\physBoundary_{D_1}$ and likewise for members of $\mathcal{E}_{D_2}$, i.e., for $\alpha = 1,2$,
\begin{equation*}
 \mathcal{E}_{D_\alpha} \equiv \left\{ E \in \mathcal{E} \colon E \subset \physBoundary_{D_\alpha} \right\}.
\end{equation*}
To ensure that each edge in $\mathcal{E}$ belongs to either the Neumann or Dirichlet boundaries, assume that $\physBoundary_{D_i} = \text{int}(\overline{\cup_{E \in \mathcal{E}_{D_\alpha}} E})$ for $\alpha = 1,2$.  We associate an edge size $h_E = h_K$ for each edge $E \in \mathcal{E}$, where $K \in \mathcal{K}$ is the element for which $E$ is the edge. With these definitions in place, we have the following lemma:

\begin{lemma}[Trace Inequalities]
  There exists two positive, dimensionless constants $\cTrace{,1}^{\B}, \cTrace{,2}^{\B} > 0$ such that
  \begin{equation}
    \sum_{E_1 \in \mathcal{E}_{D_1}} \int_{E_1} \frac{h_{E_1}^3}{\cTrace{,1}^{\B}} \left| \left( \nabla \Delta {\bf v}_h \right) \cdot {\bf n} \right|^2 \ d \physBoundary \le \frac{1}{2} a^{\B}({\bf v}_{h},{\bf v}_{h})
    \label{eqn:TI_B_1}
  \end{equation}
  \begin{equation}
   \sum_{E_2 \in \mathcal{E}_{D_2}} \int_{E_2} \frac{h_{E_2}}{\cTrace{,2}^{\B}} \left| \Delta {\bf v}_h \right|^2 \ d \physBoundary \le \frac{1}{2} a^{\B}({\bf v}_{h},{\bf v}_{h})
   \label{eqn:TI_B_2}
  \end{equation}
  for all ${\bf v}_{h} \in \mathcal{V}^{\B}_h$. Furthermore, explicit upper bounds for both $\cTrace{,1}^{\B}$ and $\cTrace{,2}^{\B}$ can be obtained through associated generalized eigenvalue problem.

  \begin{proof}
    The proof of the result for \eqref{eqn:TI_B_1} requires some additional mathematical machinery in comparison to the proof of Lemma~\ref{lemma:trace_heat_prelim}. This is because \eqref{eqn:TI_B_1} is actually a composition of two sequential relationships: first, a trace inequality, and second, an \emph{inverse inequality}. The familiar trace inequality relates the integrand on the boundary to the interior of the domain through an appropriate mesh parameter scaling, while the inverse inequality reduces the order of the derivative operator to match that of $a^{\B}( \cdot , \cdot )$. In general, inverse inequalities are used to bound high-ordered derivatives by low-ordered derivatives and hold only for finite-dimensional spaces; the ``inverse'' is in reference to the classic Poincar\'e-Friedrichs inequality \cite{harari1992}.

    To begin, we consider the kernel of the gradient of the Laplacian, which is denoted
    \begin{equation*}
      \text{ker}(\nabla \Delta) \equiv \left\{ {\bf v}_{h} \in \mathcal{V}_h^{\B} \colon \nabla \Delta {\bf v}_{h} = {\bf 0} \right\},
    \end{equation*}
    and the orthogonal compliment that comprises $\mathring{\mathcal{V}}^\B_{1,h}$. Like the approach taken in Lemma~\ref{lemma:trace_heat_prelim}, we can decompose the discrete space through the orthogonal direct sum as $\mathcal{V}^\B_h = \mathring{\mathcal{V}}^\B_{1,h} \oplus \text{ker}(\nabla \Delta)$ and we are once again afforded the same conveniences associated with trace constants.

    The trace problem of interest relates the gradient of the Laplacian on the boundary to the interior. Specifically, we are interested in the trace inequality constant characterized by the following eigenproblem: Find $({\bf u}_h, \lambda_{1,h}) \in \mathring{\mathcal{V}}^\B_{1,h} \times \R$ such that
    \begin{equation}
      \sum_{E_1 \in \mathcal{E}_{D_1}} \int_{E_1} h_{E_1} \left( \left( \nabla \Delta {\bf u}_h \right) \cdot {\bf n} \right) \cdot \left( \left( \nabla \Delta \delta {\bf u}_h \right) \cdot {\bf n} \right) \ d \physBoundary = \lambda_{1,h} \sum_{K \in \mathcal{K}} \int_K \left( \nabla \Delta {\bf u}_h \right) \cdot \left( \nabla \Delta \delta {\bf u}_h \right) \ d\physDomain
    \end{equation}
    for all $\delta {\bf u}_h \in \mathring{\mathcal{V}}^\B_{1,h}$. For this eigenproblem, we split the interior integration into a summation of element-wise integrals for much of the same reasons as those given in Section~\ref{sec:domain_enlargement} namely, to avoid the imposition of additional smoothness assumptions on the entities in $\mathring{\mathcal{V}}^\B_{1,h}$. Also note that for sake of brevity, we refrain from defining an auxiliary bilinear form here, in comparison to the proof of Lemma~\ref{lemma:trace_heat_prelim}. The min-max theorem yields that the max eigenvalue associated with this eigenproblem is given by
    \begin{equation}
      \lambda^{\B}_{1,\textup{max}} = \sup_{\substack{{\bf v}_h \in \mathring{\mathcal{V}}^{\B}_{1,h} \\ {\bf v}_h \neq {\bf 0}}} \frac{ \sum_{E_1 \in \mathcal{E}_{D_1}} \int_{E_1} h_{E_1} \left| \left( \nabla \Delta {\bf v}_h \right) \cdot {\bf n} \right|^2  \ d \physBoundary }{\sum_{K \in \mathcal{K}} \int_K \left| \nabla \Delta {\bf v}_h \right|^2 \ d\physDomain}.
      \label{eqn:RQ_B_1}
    \end{equation}
    Note that this Rayleigh quotient is well-defined since the bilinear form associated with the denominator is coercive over $\mathring{\mathcal{V}}^\B_{1,h}$.

    Finally, we employ the aforementioned inverse inequality
    \begin{equation}
      \sum_{K \in \mathcal{K}} \int_K h^2_K \left| \nabla \Delta {\bf v}_h \right|^2 \ d\physDomain \le C^\B_{\text{inv}} a^{\B}({\bf v}_{h},{\bf v}_{h})
    \end{equation}
    that holds for all ${\bf v}_h \in \mathcal{V}^\B_h$. Leveraging our edge mesh construction, i.e., $h_K \equiv h_E$, and combining the auxiliary trace inequality with this inverse inequality completes the proof of \eqref{eqn:TI_B_1} with $\cTrace{,1}^{\B} = 2 \lambda^{\B}_{1,\textup{max}} C^\B_{\text{inv}}$.

    The proof of the result for \eqref{eqn:TI_B_2} follows in an analogous manner to the first part of the proof of \eqref{eqn:TI_B_1} and to that of Lemma~\ref{lemma:trace_heat_prelim}. In this instance, we are interested in the kernel of the Laplace operator that comprises the bilinear form \eqref{eqn:aB}
    \begin{equation*}
      \text{ker}(\Delta) \equiv \left\{ {\bf v}_{h} \in \mathcal{V}_h^{\B} \colon \Delta {\bf v}_{h} = 0 \right\},
    \end{equation*}
    and the orthogonal compliment that comprises $\mathring{\mathcal{V}}^\B_{2,h}$. This permits the decomposition of $\mathcal{V}^\B_h = \mathring{\mathcal{V}}^\B_{2,h} \oplus \text{ker}(\Delta)$ by a similar reasoning to that of Lemma~\ref{lemma:trace_heat_prelim}. Accordingly, the constant $\cTrace{,2}^{\B}$ for which \eqref{eqn:TI_B_2} holds over all $\mathring{\bf v}_h \in \mathring{\mathcal{V}}^\B_{2,h}$ also holds for all ${\bf v}_h \in \mathcal{V}^\B_h$, since $\Delta {\bf v}_h = \Delta \mathring{\bf v}_h$.

    Following the steps taken in Lemma~\ref{lemma:trace_heat_prelim}, we consider a generalized eigenproblem associated with \eqref{eqn:TI_B_2}, akin to that of \eqref{eq:Gen_EP_M}, and by the min-max theorem, the associated max eigenvalue is given by a Rayleigh quotient,
    \begin{equation}
      \lambda^{\B}_{2,\textup{max}} = \sup_{\substack{{\bf v}_h \in \mathring{\mathcal{V}}^{\B}_{2,h} \\ {\bf v}_h \neq {\bf 0}}} \frac{ \sum_{E_2 \in \mathcal{E}_{D_2}} \int_{E_2} h_{E_2} \left| \Delta {\bf v}_h \right|^2 \ d \physBoundary }{a^{\B} ( {\bf v}_h, {\bf v}_h ) }.
      \label{eqn:RQ_B}
    \end{equation}
    This leads to the conclusion that \eqref{eqn:TI_B_2} is satisfied for $\cTrace{,2}^{\B} = 2 \lambda^{\B}_{2,\textup{max}}$.

  \end{proof}
  \label{lemma:TI_B}
\end{lemma}

\begin{remark}
  The coefficient of $\tfrac{1}{2}$ in both \eqref{eqn:TI_B_1} and \eqref{eqn:TI_B_2} is added to result in a cleaner final Nitsche formulation however the same arguments therein hold if the $\tfrac{1}{2}$ is removed.
\end{remark}

\begin{remark}
  In the proof of Lemma~\ref{lemma:TI_B}, we specifically stated that the trace inequality followed by an inverse inequality are applied in sequence. This subtle detail is important because the inverse inequality does not necessarily hold on the domain boundary. For instance, it is possible to construct ${\bf v}_h$ such that $\nabla \Delta {\bf v}_h$ is not finite while $\Delta {\bf v}_h < \infty$ thus violating the relationship.
\end{remark}

\begin{remark}
  Note that the kernel of the Laplacian used in the proof of Lemma~\ref{lemma:TI_B} is quite a bit more complex than that of the gradient. For this reason, we make no attempt to concisely characterize it as we did in Lemma~\ref{lemma:trace_heat_prelim}.
\end{remark}

\begin{remark}
  As stated in the proof of Lemma~\ref{lemma:TI_B}, \eqref{eqn:TI_B_1} and \eqref{eqn:TI_B_2} are satisfied for $C^{\B}_{\textup{tr},1} = 2 \lambda^{\B}_{1,\textup{max}} C^\B_{\text{inv}}$ and $C^{\B}_{\textup{tr},2} = 2 \lambda^{\B}_{\textup{max},2}$, respectively. However, similar to Remark~\ref{remark:model_EP_for_TI}, constructing the basis functions that span $\mathring{\mathcal{V}}^\B_{1,h}$ and $\mathring{\mathcal{V}}^\B_{2,h}$ is impractical. Fortunately, by the reasoning of Remark~\ref{remark:model_EP_for_TI}, we can utilize the same approach of solving the associated eigenproblems over the entire discrete spaces, $\mathcal{V}^\B_{1,h}$ and $\mathcal{V}^\B_{2,h}$ instead. Furthermore, for \eqref{eqn:TI_B_1}, we are also able to embed both the trace inequality and inverse inequality into the following eigenproblem: Find $({\bf u}_{h}, \lambda_h) \in \mathcal{V}^{\B}_h \times \mathbb{R}$ such that
  \begin{equation*}
  \sum_{E_1 \in \mathcal{E}_{D_1}} \int_{E_1} h_{E_1}^3 \left[ \left( \nabla \Delta {\bf u}_h \right) \cdot {\bf n} \right] \cdot \left[ \left( \nabla \Delta \delta {\bf u}_h \right) \cdot {\bf n} \right] \ d \physBoundary = \frac{\lambda_h}{2} a^{\B}({\bf u}_{h},\delta {\bf u}_{h})
  \end{equation*}
  for all $\delta {\bf u}_{h} \in \mathcal{V}^{\B}_h$.
\end{remark}

To successfully arrive at a Nitsche formulation for the biharmonic problem, we must specify suitable linear maps such that the generalized trace and Cauchy-Schwarz inequalities appearing in Assumption~\ref{assumption2} are satisfied. The boundary operator $\mathcal{B}^{\B} \allowbreak \colon \tilde{\mathcal{V}}^{\B} \rightarrow \left( \mathcal{Q}^{\B} \right)^*$ was implicitly defined in Lemma~\ref{lemma:greens_biharmonic} for ${\bf w} \in \tilde{\mathcal{V}}^{\B}$ and ${\bf v} \in \mathcal{V}^{\B}$. According to our abstract framework, we must extend the domain of definition of this operator to the enlarged space $\tilde{\mathcal{V}}^{\B} + \mathcal{V}^{\B}_h$. As was done for the vector-valued Poisson problem, this is accomplished by expressing the boundary operator as a summation of integrals over element edges. In particular, the following expression is well-defined for any $C^1$-continuous polynomial or rational approximation over $\mathcal{K}$:
\begin{equation}
  \begin{aligned}
\left\langle \mathcal{B}^{\B} {\bf w}, \mathcal{T}^{\B} {\bf v} \right\rangle = - \sum_{E_1 \in \mathcal{E}_{D_1}} \int_{E_1} \left[ \left( \nabla \Delta {\bf w} \right) \cdot {\bf n} \right] \cdot {\bf v} d \physBoundary + \sum_{E_2 \in \mathcal{E}_{D_2}} \int_{E_2} \left[ \left( \nabla {\bf v} \right) \cdot {\bf n} \right] \cdot \left( \Delta {\bf w} \right) d \physBoundary.
  \end{aligned}
  \label{eqn:B_BvTw}
\end{equation}

Let $\eta^{\B} \colon \text{dom}(\eta^{\B}) \subseteq \left( \mathcal{Q}^{\B} \right)^* \rightarrow \mathcal{Q}^{\B}$ be a densely defined, positive, self-adjoint linear map that is defined on the enlarged space
\begin{equation*}
  \left\{ \mathcal{B}^{\B} {\bf v} \colon {\bf v} \in \tilde{\mathcal{V}}^{\B} + \mathcal{V}^{\B}_h \right\}
\end{equation*}
and satisfies
\begin{equation}
  \begin{aligned}
    \big\langle \left. \mathcal{B}^{\B} {\bf w} , \eta^{\B} \mathcal{B}^{\B} {\bf v} \right\rangle &=
    \sum_{E_1 \in \mathcal{E}_{D_1}} \int_{E_1} \frac{h_{E_1}^3}{\cTrace{,1}^{\B}} \left[ \left( \nabla \Delta {\bf w} \right) \cdot {\bf n} \right] \cdot \left[ \left( \nabla \Delta {\bf v} \right) \cdot {\bf n} \right] d \physBoundary\\
    &+ \sum_{E_2 \in \mathcal{E}_{D_2}} \int_{E_2} \frac{h_{E_2}}{\cTrace{,2}^{\B}} \left( \Delta {\bf w} \right) \cdot \left( \Delta {\bf v} \right) d \physBoundary
  \end{aligned}
  \label{eqn:eta_B}
\end{equation}
for all ${\bf w}, {\bf v} \in \tilde{\mathcal{V}}^{\B} + \mathcal{V}^{\B}_h$. Similar to that of the vector Poisson problem, the implicit definition for $\eta^{\B}$ is inspired by the results of Lemma~\ref{lemma:TI_B}.

\begin{lemma}[Generalized Trace Inequality for the Biharmonic Problem]
  Using the definition in \eqref{eqn:eta_B}, it holds that
  \begin{equation*}
    \left\langle \mathcal{B}^{\B} {\bf v}_{h}, \eta^{\B} \mathcal{B}^{\B} {\bf v}_{h} \right\rangle \leq a^{\B}({\bf v}_{h},{\bf v}_{h})
  \end{equation*}
  for all ${\bf v}_{h} \in \mathcal{V}^{\B}_h$.
  \begin{proof}
    The proof follows immediately from Lemma~\ref{lemma:TI_B} and the definition of $\eta^{\B}$.
  \end{proof}
  \label{lemma:TI_B_gen}
\end{lemma}

Next, we define the linear map $\epsilon^{\B} \colon \textup{dom}(\epsilon^{\B}) \subseteq \left( \mathcal{Q}^{\B} \right)^* \rightarrow \mathcal{Q}^{\B}$ through the action of its inverse as
\begin{equation}
  \begin{aligned}
    \Big\langle \left( \epsilon^{\B} \right)^{-1} {\bf w}, {\bf v} \Big\rangle
    &:= \sum_{E_1 \in \mathcal{E}_{D_1}} \int_{E_1} \frac{\cPen{,1}^{\B}}{h_{E_1}^3} {\bf w} \cdot {\bf v} \ d \physBoundary\\ &+ \sum_{E_2 \in \mathcal{E}_{D_2}} \int_{E_2} \frac{\cPen{,2}^{\B}}{h_{E_2}} \left[ \left( \nabla {\bf w} \right) \cdot {\bf n} \right] \cdot \left[ \left( \nabla {\bf v} \right) \cdot {\bf n} \right] \ d \physBoundary,
  \end{aligned}
  \label{eqn:eps_B}
\end{equation}
for all ${\bf w}, {\bf v} \in \mathcal{Q}^{\B}$, where $\cPen{,1}^{\B} > \cTrace{,1}^{\B}$ and $\cPen{,2}^{\B} > \cTrace{,2}^{\B}$ are positive dimensionless constants.

\begin{lemma}[Generalized Cauchy-Schwarz Inequality for the Biharmonic Problem]
  Let $\cPen{,1}^{\B} = \CSconstant_1^2 \cTrace{,1}^{\B}$ and $\cPen{,2}^{\B} = \CSconstant_2^2 \cTrace{,2}^{\B}$, where $\CSconstant_1, \CSconstant_2 \in (1,\infty)$. Then
  \begin{equation*}
    \left| \left\langle \mathcal{B}^{\B} {\bf v}, \mathcal{T}^{\B} {\bf w} \right\rangle \right| \le \frac{1}{\CSconstant} \left\langle \mathcal{B}^{\B} {\bf v}, \eta^{\B} \mathcal{B}^{\B} {\bf v} \right\rangle^{1/2} \left\langle \left( \epsilon^{\B} \right)^{-1} \mathcal{T}^{\B} {\bf w}, \mathcal{T}^{\B} {\bf w} \right\rangle^{1/2}
  \end{equation*}
  for all ${\bf v}, {\bf w} \in \tilde{\mathcal{V}}^{\B} + \mathcal{V}^{\B}_h$, where $\CSconstant = \min(\CSconstant_1,\CSconstant_2)$.
  \begin{proof}
    We individually bound the two terms in \eqref{eqn:B_BvTw} analogously to the proof of Lemma~\ref{lemma:cs_heat} by utilizing standard continuous and discrete Cauchy-Schwarz inequalities. Thereafter, the desired result is obtained by combining these bounds with the additional bounds $1/\CSconstant_1, 1/\CSconstant_2 < 1/\CSconstant$, where $\CSconstant = \min(\CSconstant_1,\CSconstant_2)$, followed by an application of the discrete Cauchy-Schwarz inequality. Much like the proof of the generalized Cauchy-Schwarz inequality for the vector Poisson problem, the expressions for $\cPen{,1}^{\B}$ and $\cPen{,2}^{\B}$ arise from the relationships
    \begin{equation}
      1 = \frac{\cPen{,1}^{\B}}{\CSconstant_1^2 \cTrace{,1}^{\B}} = \frac{1}{\CSconstant_1^2} \frac{h^3_{E_1}}{\cTrace{,1}^{\B}} \frac{\cPen{,1}^{\B}}{h^3_{E_1}} \hspace{10pt} \text{and} \hspace{10pt} 1 = \frac{\cPen{,2}^{\B}}{\CSconstant_2^2 \cTrace{,2}^{\B}} = \frac{1}{\CSconstant_2^2} \frac{h_{E_2}}{\cTrace{,2}^{\B}} \frac{\cPen{,2}^{\B}}{h_{E_2}}
    \end{equation}
  \end{proof}
  \label{lemma:CS_B}
\end{lemma}
After these choices of linear maps have been made, the generalized trace and Cauchy-Schwarz inequalities appearing in Assumption~\ref{assumption2} are satisfied and we are ready to present Nitsche's method for the biharmonic problem.


\subsection{Nitsche's Method}

Following the abstract variational framework of Section~\ref{sec:Nitsche} and with the appropriate definitions of $\epsilon^{\B}$, $\eta^{\B}$, and $\mathcal{B}^{\B}$ in place, Nitsche's method for the biharmonic equation is posed as follows:

\vspace{1em}
\begin{mybox}[\emph{Nitsche's Method for the Biharmonic Problem}]
  \vspace{-3pt}
$$
(N^{\B}_h) \left\{ \hspace{5pt}
\parbox{4.4in}{
Given $\linFunctional^{\B} \in \left( \mathcal{V}^{\B} \right)^*$ and $\left( \textup{\bf g},\textup{\bf h} \right) \in \mathcal{Q}^{\B}$, find ${\bf u}_{h} \in \mathcal{V}^{\B}_h$ such that
\begin{equation*}
\begin{aligned}
&a_h^{\B} ({\bf u}_{h},\delta {\bf u}_{h}) = \underbrace{ \int_{\physDomain} \textup{\bf f} \cdot \delta {\bf u}_{h} d \physDomain - \int_{\physBoundary_{N_1}} \textup{\bf p} \cdot \delta {\bf u}_h d \physBoundary + \int_{\physBoundary_{N_2}} \textup{\bf q} \cdot \left[ \left( \nabla \delta {\bf u}_h \right) \cdot {\bf n} \right] d \physBoundary }_{ \displaystyle \langle \linFunctional^{\B}, \delta {\bf u}_{h} \rangle}\\
&{\color{ForestGreen} \underbrace{ + \hspace{-4pt} \sum_{E_1 \in \mathcal{E}_{D_1}} \hspace{-3pt} \int_{E_1} \hspace{-3pt} \left[ \left( \nabla \Delta \delta {\bf u}_{h} \right) \cdot {\bf n} \right] \cdot \textup{\bf g} \ d \physBoundary - \hspace{-4pt} \sum_{E_2 \in \mathcal{E}_{D_2}} \hspace{-3pt} \int_{E_2} \hspace{-3pt} \left( \Delta \delta {\bf u}_{h} \right) \cdot \textup{\bf h} \ d \physBoundary }_\text{Symmetry Terms} }\\
&{\color{Orchid} \underbrace{ + \hspace{-4pt} \sum_{E_1 \in \mathcal{E}_{D_1}} \hspace{-3pt} \int_{E_1} \hspace{-3pt} \frac{\cPen{,1}^{\B}}{h_{E_1}^3} \delta {\bf u}_{h} \cdot \textup{\bf g} \ d \physBoundary + \hspace{-4pt} \sum_{E_2 \in \mathcal{E}_{D_2}} \hspace{-3pt} \int_{E_2} \hspace{-3pt} \frac{\cPen{,2}^{\B}}{h_{E_2}} \left[ \left( \nabla \delta {\bf u}_{h} \right) \cdot {\bf n} \right] \cdot \textup{\bf h} d \physBoundary }_\text{Penalty Terms} }
\end{aligned}
\label{eqn:B_Weak_Nitsche}
\end{equation*}
for every $\delta {\bf u}_{h} \in \mathcal{V}^{\B}_h$, where $a^{\B}_h: \left( \tilde{\mathcal{V}}^{\B} + \mathcal{V}^{\B}_h \right) \times \left( \tilde{\mathcal{V}}^{\B} + \mathcal{V}^{\B}_h \right) \rightarrow \mathbb{R}$ is the bilinear form defined by
\begin{equation*}
\begin{aligned}
&a^{\B}_h ({\bf w}_{h}, {\bf v}_{h}) \equiv \underbrace{ \int_{\physDomain} \left( \Delta {\bf w}_{h} \right) \cdot \left( \Delta {\bf v}_{h} \right) d \physDomain }_{ \displaystyle a^{\B}({\bf w}_{h}, {\bf v}_{h}) } \\
& {\color{Cerulean} \underbrace{ + \hspace{-4pt} \sum_{E_1 \in \mathcal{E}_{D_1}} \hspace{-3pt} \int_{E_1} \hspace{-3pt} \left[ \left( \nabla \Delta {\bf w}_{h} \right) \cdot {\bf n} \right] \cdot {\bf v}_{h} d \physBoundary - \hspace{-4pt} \sum_{E_2 \in \mathcal{E}_{D_2}} \hspace{-3pt} \int_{E_2} \hspace{-3pt} \left( \Delta {\bf w}_{h} \right) \cdot \left[ \left( \nabla {\bf v}_h \right) \cdot {\bf n} \right] d \physBoundary }_\text{Consistency Terms} }\\
& {\color{ForestGreen} \underbrace{ + \hspace{-4pt} \sum_{E_1 \in \mathcal{E}_{D_1}} \hspace{-3pt} \int_{E_1} \hspace{-3pt} \left[ \left( \nabla \Delta {\bf v}_{h} \right) \cdot {\bf n} \right] \cdot {\bf w}_{h} d \physBoundary - \hspace{-4pt} \sum_{E_2 \in \mathcal{E}_{D_2}} \hspace{-3pt} \int_{E_2} \hspace{-3pt} \left( \Delta {\bf v}_{h} \right) \cdot \left[ \left( \nabla {\bf w}_h \right) \cdot {\bf n} \right] d \physBoundary }_\text{Symmetry Terms} }\\
& {\color{Orchid} \underbrace{ + \hspace{-4pt} \sum_{E_1 \in \mathcal{E}_{D_1}} \hspace{-3pt} \int_{E_1} \hspace{-3pt} \frac{\cPen{,1}^{\B}}{h_{E_1}^3} {\bf w}_{h} \cdot {\bf v}_{h} d \physBoundary + \hspace{-4pt} \sum_{E_2 \in \mathcal{E}_{D_2}} \hspace{-3pt} \int_{E_2} \hspace{-3pt} \frac{\cPen{,2}^{\B}}{h_{E_2}} \left[ \left( \nabla {\bf v}_h \right) \cdot {\bf n} \right] \cdot \left[ \left( \nabla {\bf w}_h \right) \cdot {\bf n} \right] d \physBoundary }_{\text{Penalty Terms}} }
\end{aligned}
\end{equation*}
for ${\bf w}_h, {\bf v}_h \in \tilde{\mathcal{V}}^{\B} + \mathcal{V}^{\B}_h$.
}
\right.
\vspace{3pt}
$$
\end{mybox}
\vspace{1em}

\noindent Since we have constructed Nitsche's method such that Assumptions~\ref{assumption1} and~\ref{assumption2} are satisfied according to Lemmas~\ref{lemma:greens_biharmonic},~\ref{lemma:TI_B_gen}, and~\ref{lemma:CS_B}, we have the following theorem stating well-posedness and an error estimate for Nitsche's method for the biharmonic equation:

\begin{theorem}[Well-Posedness and Error Estimate for the Biharmonic Problem]
  Let $\cPen{,1}^{\B} = \CSconstant_1^2 \cTrace{,1}^{\B}$ and $\cPen{,2}^{\B} = \CSconstant_2^2 \cTrace{,2}^{\B}$, where $\CSconstant_1, \CSconstant_2 \in (1,\infty)$. Then there exists a unique discrete solution ${\bf u}_{h} \in \mathcal{V}^{\B}_h$ to the Nitsche formulation of the Biharmonic Problem $(N^{\B}_h)$. Moreover, if the continuous solution ${\bf u} \in \mathcal{V}^{\B}$ to Problem $(V^{\B})$ satisfies ${\bf u} \in \tilde{\mathcal{V}}^{\B}$, then the discrete solution ${\bf u}_{h}$ satisfies the error estimate
  \begin{equation*}
\vvvertiii{{\bf u} - {\bf u}_{h}}_{\B} \leq \left( 1+ \frac{2}{1-\frac{1}{\CSconstant}} \right) \min_{{\bf v}_{h} \in \mathcal{V}^{\B}_h} \vvvertiii{{\bf u} - {\bf v}_{h}}_{\B},
\label{eqn:B_Error}
\end{equation*}
  where $\CSconstant = \min(\CSconstant_1,\CSconstant_2)$.
  \begin{proof}
    The proof follows identically to that of Theorem \ref{theorem:convergence_heat}.
  \end{proof}
\end{theorem}

\begin{remark}
The choice of penalty constants presented in this subsection is not the only stable choice.  For discretization-dependent dimensionless constants $\alpha_1 > 1$, and $\alpha_2 > 1$, one can alternatively select $\cPen{,1}^{\B} > \alpha_1 \cTrace{}^{\B}$ and $\cPen{,2}^{\B} > \alpha_2 \cTrace{}^{\B}$, where $\cTrace{}^{\B} > 0$ is a dimensionless constant such that
  \begin{equation*}
  \begin{aligned}
    \left( \sum_{E_1 \in \mathcal{E}_{D_1}} \int_{E_1} \frac{h_{E_1}^3}{\alpha_1} \left| \left( \nabla \Delta {\bf v}_h \right) \cdot {\bf n} \right|^2 d \physBoundary + \sum_{E_2 \in \mathcal{E}_{D_2}} \int_{E_2} \frac{h_{E_2}}{\alpha_3 } \left| \Delta {\bf v}_h \right|^2 d \physBoundary \right) &\le \cTrace{}^{\B} a^{\B}({\bf v}_{h},{\bf v}_{h})
  \end{aligned}
  \end{equation*}
  for all ${\bf v}_{h} \in \mathcal{V}^{\B}_h$.  The advantage of this approach is that only one trace constant, namely, $\cTrace{}^{\B}$, must be estimated.  The disadvantage of this approach is that $\alpha_1$ and $\alpha_2$, which control the relative weightings of the displacement boundary condition along $\Gamma_{D_1}$ and the rotation boundary condition along $\Gamma_{D_2}$, respectively, must be specified.
\label{remark:eig_B}
\end{remark}

This concludes our derivation of Nitsche's method for the biharmonic problem. Among other things, the process has demonstrated that our abstract framework is robust with respect to the order of the PDE. Although there are more complex steps involved in arriving at Nitsche formulation for this problem, the process is entirely procedural and follows the same steps as for the Poisson problem in Section~\ref{sec:Poisson}. Unlike that for the Poisson problem, the resulting formulation here contains two sets of symmetry, consistency, and penalty terms to handle the two sets of boundary conditions that arise in the biharmonic problem.


\section{Nitsche's Method for the linearized Kirchhoff-Love Plate}
\label{sec:KL_plate}

We finish by considering Nitsche's method for the linearized Kirchhoff-Love plate. This can be obtained through a slight modification of the biharmonic problem. While the two problems are closely related, the plate problem is interesting because it utilizes a set of physical boundary conditions that are not naturally admissible to our framework. We therefore show how to manipulate these boundary conditions in such a way that they are amenable to our methodology.

A general plate model employs a transverse displacement variable, denoted $u_z$, as well as a rotational degree of freedom, denoted $\theta$. The Kirchhoff-Love kinematic assumption is that a straight line normal to the mid-surface after deformation remains straight, unstretched, and normal to the deformed surface. Physically, this asserts that the resulting displacement field is free of transverse shear strain, which introduces a constraint between the midsurface rotational variable and the transverse displacement degree of freedom, namely, $\theta(w_z) = - \nabla w_z$. Along the boundary, this rotation can be decomposed further into the \textbf{\emph{normal rotation}}, $\theta_n(w_z) = - \nabla w_z \cdot {\bf n}$, where ${\bf n}$ is the outward-facing normal on $\physBoundary$, and the \textbf{\emph{twisting rotation}}, $\theta_t(w_z) = - \nabla w_z \cdot {\bf t}$, where ${\bf t}$ is the postively-oriented tangent along $\physBoundary$.

To formulate Nitsche's method for the Kirchhoff-Love plate, we set $\sDim=2$, $\uDim=1$, $\mathcal{V}^{\KLP} \equiv \mathcal{V}^{\B}$, and $\mathcal{Q}^{\KLP} \equiv \mathcal{Q}^{\B}$, where these spaces are defined in Section~\ref{sec:biharmonic}, with superscript $\KLP$ denoting quantities associated with the Kirchhoff-Love plate problem. The boundary, $\physBoundary$, is decomposed into parts identically to that of the biharmonic. For reasons that will be apparent later, we additionally define the set $\cornerSet{} \subset \physBoundary$ as the set of ``corners''.  We further decompose this set into $\cornerSet{D} \equiv \cornerSet{} \cap \overline{\physBoundary_{D_1}}$ and $\cornerSet{N} \equiv \cornerSet{} \cap \physBoundary_{N_1}$ and note that, by construction, $\cornerSet{} = \cornerSet{D} \cup \cornerSet{N}$ and $\cornerSet{D} \cap \cornerSet{N} = \emptyset$.  We denote corners as $C \in \cornerSet{}$.

We utilize stress and strain measures in the plate model as a proxy for the resulting displacement field since they are more convenient for formulating the variational problem. The \textbf{\emph{bending strain}} for the Kirchhoff-Love plate is defined as $\bm{\bendStrain}(w_z) \equiv -\nabla^s \nabla w_z$, where $\nabla^s$ is the symmetric part of the gradient. The \textbf{\emph{bending stress}}, ${\bf B}$, for an isotropic material can be expressed via Hooke's law as ${\bf B} = D [\nu \textup{tr}(\bm{\bendStrain}) {\bf I} + (1-\nu) \bm{\bendStrain} ]$, where $D = \thickness^3 E / [12 (1-\nu^2)]$ is the isotropic plate rigidity, $E$ is Young's modulus, $\nu$ is Poisson's ratio, $\thickness$ is the plate thickness, ${\bf I}$ is the identity tensor, and $\textup{tr}(\cdot)$ denotes the tensor trace. The internal energy associated with the linearized Kirchhoff-Love plate is given by the strain energy due to plate bending, namely
\begin{equation}
    E_{\textup{in}}^{\KLP}(v_z) = \frac{1}{2} \int_{\physDomain} {\bf B}(v_z) : \bm{\bendStrain}(v_z) \ d \physDomain.
\end{equation}

Given this, we readily define an associated bilinear form $a^{\KLP}(\cdot,\cdot) \colon \mathcal{V}^{\KLP} \times \mathcal{V}^{\KLP} \rightarrow \mathbb{R}$ via
\begin{equation*}
  a^{\KLP}(w_z,v_z) \equiv \int_{\physDomain} {\bf B}(w_z) : \bm{\bendStrain}(v_z) \ d \physDomain
  \label{eqn:P_int}
\end{equation*}
for all $w_z, v_z \in \mathcal{V}^{\KLP}$.

Let $\applied{\forcing}_z \in L^2(\physDomain)$ be the applied transverse loading, $\applied{u}_z \in H^{3/2}(\physBoundary_{D_1})$ the prescribed transverse displacement, and $\applied{\theta}_n \in H^{1/2}(\physBoundary_{D_2})$ the prescribed normal rotation. Accordingly, we define the trace operator $\mathcal{T}^{\KLP} \colon \mathcal{V}^{\KLP} \rightarrow \mathcal{Q}^{\KLP}$ via its action on the transverse displacement $v_z \in \mathcal{V}^{\KLP}$, i.e., $\mathcal{T}^{\KLP} v_z = \left( v_z \big|_{\physBoundary_{D_1}}, \midsurfRot_n(v_z) \big|_{\physBoundary_{D_2}} \right)$. Given $\left( \applied{u}_z, \applied{\midsurfRot}_n \right) \in \mathcal{Q}^{\KLP}$, we define
\begin{equation*}
  \mathcal{V}^{\KLP}_{\applied{u}_z,\applied{\midsurfRot}_n} \equiv \left\{  v_z \in \mathcal{V}^{\KLP} \colon \mathcal{T}^{\KLP} v_z = \left( \applied{u}_z, \applied{\midsurfRot}_n \right) \right\}
\end{equation*}
as the trial space of displacement fields satisfying the prescribed Dirichlet boundary conditions while $\mathcal{V}^{\KLP}_{0,0}$ is the homogeneous counterpart that comprises the test space.

The Neumann boundary conditions resemble that of the biharmonic, namely, an applied transverse shearing, $\applied{\tau}_z \in L^2(\physBoundary_{N_1})$, and an applied moment, $\applied{\bf \moment} = \applied{\moment}_{nn} {\bf n} + \applied{\moment}_{nt} {\bf t}$ for a square-integrable bending moment, $\applied{\moment}_{nn}$, and twisting moment, $\applied{\moment}_{nt}$. Thus, it appears that the plate can accommodate one additional Neumann condition in comparison to the biharmonic, giving rise to the following form governing the external work associated with the displacement $v_z$ for the Kirchhoff-Love plate:
\begin{equation}
\begin{aligned}
E_{\textup{ext}}^{\KLP}(v_z) \equiv -\int_{\physDomain} \applied{\forcing}_z v_z \ d \physDomain - \int_{\physBoundary_{N_1}} \applied{\tau}_z v_z \ d \physBoundary - \int_{\physBoundary_{N_1}} \applied{\moment}_{nt} \midsurfRot_t(v_z) \ d \physBoundary - \int_{\physBoundary_{N_2}} \applied{\moment}_{nn} \midsurfRot_n(v_z) \ d \physBoundary.
\end{aligned}
\label{eqn:E_ext_KLP}
\end{equation}

Given the internal energy and external work, the minimization problem of interest is simply
$$
(M^{\KLP}) \left\{ \hspace{5pt}
\parbox{4.35in}{
\noindent Find $u_z \in \mathcal{V}^{\KLP}_{\applied{u}_z,\applied{\midsurfRot}_n}$ that minimizes the total energy
\begin{eqnarray*}
E^{\KLP}_{\textup{total}}(u_z) = E^{\KLP}_{\textup{int}}(u_z) + E^{\KLP}_{\textup{ext}}(u_z).
\end{eqnarray*}
}
\right.
$$

In \eqref{eqn:E_ext_KLP}, there are three boundary integrals corresponding to an applied shear or moment and their corresponding energetically conjugate displacement or rotation induced by such a force. However, our selection of $\mathcal{Q}^{\KLP}$ provides us with direct control of only two of these: the boundary displacement and normal rotation. This is problematic because in order to apply our abstract framework to the plate, we need the dual forces that correspond to enforcing these boundary conditions. In other words, this definition of $E_{\textup{ext}}^{\KLP}(v_z)$ does not admit a linear functional that is admissible to our framework because Assumption~\ref{assumption1} does not hold. The problem is that of the three proposed boundary conditions, one can be expressed as a linear combination of the other two as first observed by Kirchhoff in 1850 \cite{Kirchhoff1850}, with the corner force correction later presented by Lamb in 1889 \cite{lamb1889}.

To address this limitation, we now proceed to arrive at a linear functional that is admissible to our framework by utilizing the definition of the twisting rotation and its relation to the tangential derivative. Through integration-by-parts note that
\begin{equation}
 \int_{\physBoundary_{N_1}} \applied{\moment}_{nt} \midsurfRot_t(v_z) \ d \physBoundary = \int_{\physBoundary_{N_1}} v_z \frac{\partial \applied{\moment}_{nt}}{\partial t} \ d {\physBoundary} + \sum_{C \in \cornerSet{N}} \llbracket \applied{\moment}_{nt} \rrbracket v_z \Big|_{C}
 \label{eqn:Ersatz_IBP_KLP}
\end{equation}
for any $v_z: \physBoundary_{N_1} \rightarrow \mathbb{R}$ with $v_z \big|_{\partial \physBoundary_{N_1}} = 0$ where $\frac{\partial}{\partial t} \applied{\moment}_{nt} = \nabla \applied{\moment}_{nt} \cdot {\bf t}$. We are then able define the \textbf{\emph{ersatz traction}} via
\begin{equation*}
  \applied{\ersatz}_z \equiv \applied{\tau}_z + \frac{\partial \applied{\moment}_{nt}}{\partial t},
\end{equation*}
by combining the new boundary integral in \eqref{eqn:Ersatz_IBP_KLP} with the other boundary integral defined over $\physBoundary_{N_1}$. Furthermore, we define the \textbf{\emph{corner forces}} $\llbracket \applied{\moment}_{nt} \rrbracket$ via the jump operator
\begin{equation*}
  \llbracket \applied{\moment}_{nt}(\textbf{x}) \rrbracket \equiv \lim_{\epsilon \rightarrow 0} \left[ \applied{\moment}_{nt}(\textbf{x} + \epsilon \bdyTangent) - \applied{\moment}_{nt}(\textbf{x} - \epsilon \bdyTangent) \right].
\label{eqn:jumpDef}
\end{equation*}
Unlike the transverse shearing and twisting moment, the ersatz traction and corner forces are energetically conjugate to the boundary displacement and, hence, are the natural shears to employ in our derivation of Nitsche's method for the Kirchhoff-Love plate. We assume that $\applied{\ersatz}_z \in L^2(\physBoundary_{N_1})$ and $\left\{ \llbracket \applied{\moment}_{nt} \rrbracket \big|_{C} \right\}_{C \in \cornerSet{N}} \in \mathbb{R}^{\#\cornerSet{N}}$, respectively. Finally, this allows us to define the linear functional $\linFunctional^{\KLP} \in \left( \mathcal{V}^{\KLP} \right)^*$ that governs the external work done on the Kirchhoff-Love plate via
\begin{equation*}
\left\langle \linFunctional^{\KLP}, v_z \right\rangle \equiv \int_{\physDomain} \applied{\forcing}_z v_z \ d \physDomain + \int_{\physBoundary_{N_1}} \applied{\ersatz}_z v_z \ d \physBoundary + \sum_{C \in \cornerSet{N}} \llbracket \applied{\moment}_{nt} \rrbracket v_z \Big|_C + \int_{\physBoundary_{N_2}} \applied{\moment}_{nn} \midsurfRot_n(v_z) \ d \physBoundary \nonumber
\end{equation*}
for all $v_z \in \mathcal{V}^{\KLP}$ that is admissible to our framework. Comparing $\linFunctional^{\KLP}$ to $E_{\textup{ext}}^{\KLP}(u_z)$, we see that the ersatz traction and corner forces are indeed the Lagrange multiplier fields associated with enforcing the displacement boundary condition.

This leads to the following variational problem governing the Kirchhoff-Love plate:
$$
(V^{\KLP}) \left\{ \hspace{5pt}
\parbox{4.35in}{
\noindent Find $u_z \in \mathcal{V}^{\KLP}_{\applied{u}_z,\applied{\theta}_n}$ such that \vspace{-5pt}
\begin{eqnarray*}
  a^{\KLP}(u_z,\delta u_z) = \left\langle \linFunctional^{\KLP}, \delta u_z \right\rangle
  \label{eqn:Plate_Weak}
\end{eqnarray*}
\noindent \vspace{-5pt} for every $\delta u_z \in \mathcal{V}^{\KLP}_{0,0}$.
}
\right.
$$

Through similar procedures to what was done in earlier sections, we arrive at the following generalized Green's identity for the linearized Kirchhoff-Love plate. For simplicity in this instance, let $\tilde{\mathcal{V}}^{\KLP} \equiv H^4(\Omega) \subset \mathcal{V}^{\KLP}$. For $w_z \in \tilde{\mathcal{V}}^{\KLP}$ and $v_z \in \mathcal{V}^{\KLP}$, we have that
  \begin{equation}
    \begin{aligned}
      &a^{\KLP}(w_z,v_z)\\
      &= \underbrace{ \int_{\physBoundary_{D_2}} \moment_{nn}(w_z) \theta_n(v_z) \ d \physBoundary + \sum_{C \in \cornerSet{D}} \llbracket \moment_{nt}(w_z) \rrbracket v_z \Big|_{C} + \int_{\physBoundary_{D_1}} \ersatz_z(w_z) v_z \ d \physBoundary }_{ \displaystyle \langle \mathcal{B}^{\KLP} w_z, \mathcal{T}^{\KLP} v_z \rangle }\\
      & \underbrace{ + \hspace{-3pt} \int_{\physBoundary_{N_2}} \hspace{-10pt} \moment_{nn}(w_z) \theta_n(v_z) d \physBoundary \hspace{-1pt} + \hspace{-5pt} \sum_{C \in \cornerSet{N}} \hspace{-3pt} \llbracket \moment_{nt}(w_z) \rrbracket v_z \Big|_{C} \hspace{-3pt} + \hspace{-3pt} \int_{\physBoundary_{N_1}} \hspace{-10pt} \ersatz_z(w_z) v_z d \physBoundary \hspace{-2pt} - \hspace{-4pt} \int_\physDomain \left( \nabla \cdot \nabla \cdot {\bf B}(w_z) \right) v_z d \physDomain}_{\displaystyle \langle \mathcal{L}^{\KLP} w_z, v_z \rangle},
    \end{aligned}
    \label{eqn:Greens_ID_KLP}
  \end{equation}
  where $\moment_{nn}(w_z) \equiv {\bf n} \cdot {\bf B}(w_z) \cdot {\bf n}$ is the \textbf{bending moment}, $\moment_{nt}(w_z) \equiv {\bf n} \cdot {\bf B}(w_z) \cdot {\bf t}$ is the \textbf{twisting moment}, and
  \begin{equation*}
    \ersatz_z(w_z) \equiv \left( \nabla \cdot {\bf B}(w_z) \right) \cdot {\bf n} + \frac{\partial \moment_{nt}(w_z)}{\partial t}
    \label{eqn:KL_Plate_ersatz}
  \end{equation*}
  is the \textbf{ersatz force}.
Given the similarities between the biharmonic problem discussed in Section~\ref{sec:biharmonic} and this section, we are able to utilize the same mesh definitions. However, since the linearized Kirchhoff-Love plate requires the use of corner forces, we accordingly require the corresponding mesh entities. To this end, we associate each $C \in \cornerSet{}$ with an element $K \in \mathcal{K}$, and we define $h_C = h_K$.

From the generalized Green's identity \eqref{eqn:Greens_ID_KLP}, we are able to define $\epsilon^{\KLP}$ and $\eta^{\KLP}$ analogously to previous instances. However, a generalized Cauchy-Schwarz inequality for the linearized Kirchhoff-Love plate differs from that of the biharmonic problem presented in Section~\ref{sec:biharmonic} due to the corner force present in $\mathcal{B}^{\KLP}$. Specifically, for Assumption~\ref{assumption2}.3 to hold, we require an additional penalty term in the construction of $\epsilon^{\KLP}$ on the corner displacement that will appear in the final Nitsche formulation.

At this point, we are able pose the Nitsche formulation as follows:

\vspace{1em}
\begin{mybox}[\emph{Nitsche's Method for the Kirchhoff-Love Plate}]
  \vspace{-3pt}
$$
(N^{\KLP}_h) \left\{ \hspace{5pt}
\parbox{4.39in}{
Given $\linFunctional^{\KLP} \in \left( \mathcal{V}^{\KLP} \right)^*$ and $\left( \applied{u}_z,\applied{\midsurfRot}_n \right) \in \mathcal{Q}^{\KLP}$, find $u_{z,h} \in \mathcal{V}^{\KLP}_h$ such that
\begin{equation*}
\begin{aligned}
&a_h^{\KLP} (u_{z,h},\delta u_{z,h}) = \\
& \underbrace{ \int_{\physDomain} \hspace{-1pt} \applied{\forcing}_z \delta u_{z,h} \ d \physDomain + \int_{\physBoundary_{N_1}} \hspace{-11pt} \applied{\ersatz}_z \delta u_{z,h} \ d \physBoundary + \sum_{C \in \cornerSet{N}} \hspace{-3pt} \llbracket \applied{\moment}_{nt} \rrbracket \delta u_{z,h} \Big|_C + \int_{\physBoundary_{N_2}} \hspace{-11pt} \applied{\moment}_{nn} \midsurfRot_n(\delta u_{z,h}) \ d \physBoundary }_{ \displaystyle \langle \linFunctional^{\KLP}, \delta u_{z,h} \rangle}\\
&{\color{ForestGreen} \underbrace{ - \hspace{-8pt} \sum_{E_1 \in \mathcal{E}_{D_1}} \hspace{-5pt} \int_{E_1} \hspace{-5pt} \ersatz_z(\delta u_{z,h}) \applied{u}_z  d \physBoundary - \hspace{-5pt} \sum_{C \in \cornerSet{D}} \hspace{-3pt} \llbracket \moment_{nt}(\delta u_{z,h}) \rrbracket \applied{u}_z \Big|_{C} - \hspace{-8pt} \sum_{E_2 \in \mathcal{E}_{D_2}} \hspace{-5pt} \int_{E_2} \hspace{-5pt} \moment_{nn}(\delta u_{z,h}) \applied{\theta}_n  d \physBoundary  }_\text{Symmetry Terms} }\\
& \clipbox{-2 0 395 0}{${\color{Orchid} \underbrace{ + \thickness^3 E \left( \sum_{E_1 \in \mathcal{E}_{D_1}} \hspace{-5pt} \int_{E_1} \hspace{-3pt} \frac{\cPen{,1}^{\KLP}}{h_{E_1}^3} \delta u_{z,h} \applied{u}_z  d \physBoundary + \sum_{C \in \cornerSet{D}} \hspace{-3pt} \frac{\cPen{,2}^{\KLP}}{h_C^2} \delta u_{z,h} \applied{u}_z \Big|_C \right. \hspace{40em} }}$} \\
&\phantom{\hspace{27pt}} \clipbox{10 0 -2 0}{$ {\color{Orchid} \underbrace{ \hspace{1em} \left. + \sum_{E_2 \in \mathcal{E}_{D_2}} \hspace{-5pt} \int_{E_2} \hspace{-3pt} \frac{\cPen{,3}^{\KLP}}{h_{E_2}} \theta_n(\delta u_{z,h}) \applied{\theta}_n d \physBoundary \right) }_{\text{Penalty Terms}} }$}
\end{aligned}
\label{eqn:KL_Plate_Weak_Nitsche}
\end{equation*}
for every $\delta u_{z,h} \in \mathcal{V}^{\KLP}_h$, where $a^{\KLP}_h: \left( \tilde{\mathcal{V}}^{\KLP} + \mathcal{V}^{\KLP}_h \right) \times \left( \tilde{\mathcal{V}}^{\KLP} + \mathcal{V}^{\KLP}_h \right) \rightarrow \mathbb{R}$ is the bilinear form defined by
\begin{equation*}
\begin{aligned}
&a^{\KLP}_h(w_{z,h}, v_{z,h}) \equiv\\
& \underbrace{ \int_{\physDomain} {\bf B}(w_{z,h}) : \bm{\bendStrain}(v_{z,h}) \ d \physDomain  \hspace{-3pt} }_{ \displaystyle a^{\KLP}(w_{z,h},v_{z,h}) } \clipbox{-2 0 403 0}{ $ {\color{Cerulean} \underbrace{ - \hspace{-7pt} \sum_{E_1 \in \mathcal{E}_{D_1}} \hspace{-5pt} \int_{E_1} \hspace{-5pt} \ersatz_z(w_{z,h}) v_{z,h} d \physBoundary \hspace{-1pt} - \hspace{-7pt} \sum_{C \in \cornerSet{D}}  \hspace{-3pt} \llbracket \moment_{nt}(w_{z,h}) \rrbracket v_{z,h} \Big|_{C} \hspace{40em} } } $ }\\
 &\phantom{\hspace{27pt}}  \clipbox{10 0 -2 0}{$ {\color{Cerulean} \underbrace{ \hspace{10pt} - \hspace{-7pt} \sum_{E_2 \in \mathcal{E}_{D_2}} \hspace{-5pt} \int_{E_2} \hspace{-5pt} \moment_{nn}(w_{z,h}) \theta_n(v_{z,h}) d \physBoundary }_\text{Consistency Terms} }$ } \hspace{-5pt} \clipbox{-2 0 395 0}{ $ {\color{ForestGreen} \underbrace{ - \hspace{-5pt} \sum_{E_1 \in \mathcal{E}_{D_1}} \hspace{-5pt} \int_{E_1} \hspace{-5pt} \ersatz_z(v_{z,h}) w_{z,h} \ d \physBoundary \hspace{40em} } } $ } \\
&\phantom{\hspace{27pt}} \clipbox{10 0 -2 0}{ $ {\color{ForestGreen} \underbrace{ \hspace{10pt} - \hspace{-5pt} \sum_{C \in \cornerSet{D}}  \hspace{-3pt} \llbracket \moment_{nt}(v_{z,h}) \rrbracket w_{z,h} \Big|_{C} - \hspace{-5pt}  \sum_{E_2 \in \mathcal{E}_{D_2}} \hspace{-5pt} \int_{E_2} \hspace{-5pt} \moment_{nn}(v_{z,h}) \theta_n(w_{z,h}) \ d \physBoundary }_\text{Symmetry Terms} } $ }\\
& \clipbox{-2 0 395 0}{${\color{Orchid} \underbrace{ + \thickness^3 E \left(\sum_{E_1 \in \mathcal{E}_{D_1}} \hspace{-5pt} \int_{E_1} \hspace{-5pt} \frac{\cPen{,1}^{\KLP}}{h_{E_1}^3} v_{z,h} w_{z,h} \ d \physBoundary + \hspace{-5pt} \sum_{C \in \cornerSet{D}} \frac{\cPen{,2}^{\KLP}}{h_C^2} v_{z,h} w_{z,h} \Big|_C \right. \hspace{40em} }}$} \\
&\phantom{\hspace{27pt}} \clipbox{10 0 -2 0}{$ {\color{Orchid} \underbrace{ \hspace{1em} \left. + \hspace{-5pt} \sum_{E_2 \in \mathcal{E}_{D_2}} \hspace{-5pt} \int_{E_2} \hspace{-5pt} \frac{\cPen{,3}^{\KLP}}{h_{E_2}} \theta_n(v_{z,h}) \theta_n(w_{z,h}) \ d \physBoundary \right) }_{\text{Penalty Terms}} }$}.
\end{aligned}
\end{equation*}
for $w_{z,h}, v_{z,h} \in \tilde{\mathcal{V}}^{\KLP} + \mathcal{V}^{\KLP}_h$.
}
\right.
\vspace{3pt}
$$
\end{mybox}
\vspace{1em}

Since this formulation has been arrived at through our abstract framework, it affords the same properties as those for the biharmonic and Poisson problems, namely, well-posedness and optimal error estimates. For more detail and rigor associated with this formulation, the reader is referred to \cite{Benzaken2020}. Although that paper analyzes the Kirchhoff-Love shell, the plate counterpart is easily obtained by removing all curvature terms and considering only the out-of-plane degree of freedom. Alternatively, \cite{harari2012,gustafsson2020nitsches} derive and analyze the same Nitsche formulation, although they do not use our abstract framework.

A similar plate problem was thoroughly studied in the context of Chladni figures in \cite{gander2012chladni}. Therein, the authors emphasized the importance of the correctness of the plate formulation, and specifically, the free boundary conditions and corner forces. The problem was analyzed and solved using both spectral and finite difference discretizations. The Nitsche formulation presented in this section for the linearized Kirchhoff-Love plate offers another means for solving this problem with finite element discretizations, allowing for more general geometric configurations and higher order elements while preserving the correctness of the formulation.

In summary, we have derived a Nitsche formulation for the linearized Kirchhoff-Love plate utilizing our abstract framework. Although there is a striking resemblence of this problem to the biharmonic problem discussed in Section~\ref{sec:biharmonic}, there is a fundamental difference due to what appears as an additional physical boundary condition. This is at odds with the fact that mathematically we are still only able to enforce two conditions. To resolve this issue, we have derived the so-called ersatz forces (also known as modified shear forces), which combine the boundary displacement, or shear, with the twisting rotation, or moment. These erstatz forces imbue the operators $\mathcal{L}^{\KLP}$ and $\mathcal{B}^{\KLP}$ with the correct duality properties to satisfy Assumption~\ref{assumption1} after which posing Nitsche's method follows an identical procedure to that of the previous sections. While the ersatz forces are well-known, the fact that they follow naturally from the abstract framework shows the power of using the framework. For the shell problem discussed in \cite{Benzaken2020}, this same approach even led to the discovery of mistakes in the existing literature.


\section{Conclusion}
\label{sec:conclusion}

Learning Nitsche's method is an intimidating task for those who have not previously been exposed to this kind of formulation. The primary goal of this paper has been to provide intution for Nitsche's method as a whole, and to provide a tutorial for constructing a Nitsche formulation via the framework presented in \cite{Benzaken2020}. As we have seen, each of the extra terms required for Nitsche's method can be understood through fairly simple means. First, the penalty terms are used to weakly enforce the Dirichlet boundary conditions. Next, the consistency terms restore the variational consistency that is generally lost with the release of strongly-enforced of boundary conditions. Finally, the symmetry terms symmetrize the associated bilinear form by incorporating the residual-based, symmetric counterpart to the consistency terms.

Even with this intuition, arriving at a Nitsche formulation remains a non-trivial task. However, for variational problems, we have shown how to apply the abstract framework in a procedural fashion to arrive at a formulation that is both provably stable and convergent. The abstract framework relies on two main assumptions: (i) the existence of a generalized Green's identity that relates strong-form differential operators to the associated bilinear form, and (ii) the availability of generalized trace and Cauchy-Schwarz inequalities. Establishing results to satisfy these assumptions is the primary task required for arriving at Nitsche formulation. Once established, the framework conveniently provides a practical way of estimating the penalty constant.

Through a didactic progression of model problems with increasing complexity, we have shown how to apply the framework and rigorously arrive at the corresponding Nitsche formulations. As illustrated in the last example for the linearized Kirchhoff-Love plate, this rigor can help identify potential hurdles and devise solutions at an early stage of the construction. While the results are not new by themselves, the derivation through the abstract framework is new, and the hope is that this approach will facilitate the readers in discovering new results of their own.

\bibliographystyle{siamplain}
\bibliography{tutorial}

\end{document}